\documentclass[12pt]{iopart}  
\usepackage{amsmath}
\usepackage{amsthm}
\usepackage{amsfonts,amssymb}
\usepackage{psfig}
\usepackage{harvard}

\usepackage{graphicx,epic,eepic}   
\newtheorem{thm}{Theorem}{\bf}{\it}
\newtheorem{lemma}[thm]{Lemma}

\newtheorem{defn}[thm]{Definition}

\newtheorem*{mainexample}{Main Example}{\bf}{\it}
\def\Re{\mathrm{Re}\,}
\def\Im{\mathrm{Im}\,}
\def\Ei{\mathrm{Ei}\,}
\def\ExpEi{\mathrm{ExpEi}\,}

\def\std{\mathrm{std}\,}
\def\esssup{\mathrm{ess}\, \sup}
\renewcommand{\theequation}{\thesection.\arabic{equation}}
\newcommand{\field}[1]{\mathbb{#1}} 

\begin{document}

  \title[Renormalization for Siegel disks and Measurable Riemann Mapping Theorem]{Cylinder renormalization for Siegel disks and a constructive Measurable Riemann Mapping Theorem}

\author{Denis G. Gaidashev}
\address{Department of Mathematics, University of Toronto,
Toronto, Ontario, Canada M5S 3G3.}
\ead{gaidash@math.toronto.edu}

\begin{abstract}

The boundary of the Siegel disk of a quadratic polynomial  with an irrationally indifferent  fixed point with the golden mean rotation number has been observed to be self-similar. The geometry of this self-similarity is universal for a large class of holomorphic maps. A renormalization  explanation of this universality has been proposed in the literature. However, one of the ingredients of this explanation, the hyperbolicity of renormalization, has not been proved yet. 

The present work considers a cylinder renormalization - a novel type of renormalization for holomorphic maps with a Siegel disk which is better suited for a hyperbolicity proof. A key element of a cylinder renormalization of a holomorphic map is a conformal isomorphism of a dynamical quotient of a subset of $\field{C}$ to a bi-infinite cylinder $\field{C} /  \field{Z}$. A construction of this conformal isomorphism is an implicit procedure which can be performed using the Measurable Riemann Mapping Theorem. 

We present a constructive proof of the Measurable Riemann Mapping Theorem, and obtain rigorous bounds on a numerical approximation of the desired conformal isomorphism. Such control of the uniformizing  conformal coordinate is of  key importance for a rigorous computer-assisted study of cylinder renormalization.

\end{abstract}

\maketitle

\setcounter{page}{1}
\medskip\section{Introduction}\label{intro}
\setcounter{equation}{0}

Over the years, one particular mathematical tool -- renormalization -- has become very instrumental in problems involving universality in a class of maps or flows. 

It is possible to give a very concise definition of renormalization in the dynamical setting: A  renormalization of a dynamical system $X$ acting on a topological space $\mathcal{T}$ is a rescaled first return map of a subset of $\mathcal{T}$.

To prove universality, one usually introduces a renormalization operator on a certain functional space, and proceeds to show that a) this operator has a fixed point, and b) the operator is hyperbolic at this fixed point. 

This approach has led to dramatic successes in the rigorous explanation of the universality in several areas of complex and real dynamics, most notably, universality of unimodal maps through the works of Sullivan  \citeyear{Sul1,Sul2}, McMullen  \citeyear{McM2,McM1}, Lyubich  \citeyear{Lyub1,Lyub2} and others, universality of critical circle maps: de Faria \citeyear{dF1,dF2}, de Faria and de Melo \citeyear{dFdM1,dFdM2}, Yampolsky \citeyear{Ya1,Ya,Ya3} and others, and  universality of ``critical'' Hamiltonian flows: Koch \citeyear{Koch2,Koch,Koch3} and others. 

However, the renormalization picture for one of the most interesting objects in one-dimensional complex dynamics -- Siegel disks -- is not complete yet. Recall that a Siegel disk is the maximal linearization domain of a holomorphic map with a linearizable irrationally indifferent  fixed point.

A classical theorem of Siegel implies that for the quadratic polynomial 
\begin{equation}\label{f_theta}
f_{\theta^*}=e^{2 \pi i \theta^*} z (1-0.5 z)  
\end{equation}
with $\theta^*=(\sqrt{5}-1)/2$ (the inverse golden mean) the fixed point at the origin is linearizable. The Siegel disk ${\rm \Delta}_{\theta^*}$ of $f_{\theta^*}$ is a quasidisk whose boundary contains the critical point: $1 \in \partial {\rm \Delta}_{\theta^*}$ (cf e.g. \cite{Dou}). It has been observed numerically by \citeasnoun{MaNau} that the boundary of the Siegel disk of $f_{\theta^*}$ is asymptotically self-similar. They found that the Fibonacci iterates of the critical point, $f^{q_n}_{\theta^*}(1)$, $n=1,2,3,5,8, \ldots$, approach the critical point asymptotically along two straight lines separated by an angle $2 \alpha \approx 107.3$ degrees. At the same time, the distance between successive iterates decreases at a geometric rate:
\begin{equation}\label{Siegel_sim_1}
\lim_{n \rightarrow \infty} {f^{q_{n+1}}(1)-1 \over f^{q_n}(1)-1}=\lambda e^{(-1)^{n+1} 2 i \alpha},
\end{equation}
with $\lambda \approx -0.7419$. Thus, the geometry on the boundary of the Siegel disk repeats itself around the critical point after a reflection with respect to a certain line passing through point $1$, and a rescaling. The self-similar geometry is transported by the dynamics to other points on $\partial \Delta_{\theta^*}$.

Manton and Nauenburg also conjectured that the parameters $\alpha$ and $\lambda$, as well, as the boundary curve of the Siegel disk at the critical point and its preimages, are universal for a large class of maps with a fixed point of multiplier $ e^{2 \pi i \theta^*}$ at the origin. \citeasnoun{Wi} considered a renormalization operator defined on pairs of such commuting holomorphic maps, and conjectured the existence and hyperbolicity of a fixed point in a space of such pairs. He also obtained numerical evidence that such a fixed point does exist. 

 In the 90's there were several results that used renormalization to explain this self-similar structure. \citeasnoun{St} was able to construct a computer-assisted proof of the existence of a renormalization fixed point for holomorphic maps with a golden mean Siegel disk. More recently, \citeasnoun{McM2} proved the asymptotic self-similarity of the Siegel disk of the golden mean quadratic polynomial $f_{\theta^*}$. He has demonstrated that successive renormalizations of the golden mean quadratic converge to a fixed point. However, the renormalization for Siegel disks based on commuting pairs has not yet led to a proof of the hyperbolicity of the renormalization operator. 

The main difficulty with an attempt to prove hyperbolicity through the renormalization of commuting pairs is that the space of commuting pairs does not possess the natural structure of a Banach manifold in which renormalization is an analytic operator. Such pairs constitute a topological space whose properties are far from simple.

\citeasnoun{Ya} introduced a {\it cylinder renormalization} of analytic maps with a fixed point. The novelty of this renormalization is that it is formulated in terms of an analytic operator on a certain {\it Banach space} of analytic functions. Such an operator is more suitable for a proof of hyperbolicity. The present work fills out several steps in a computer-assisted implementation of such a proof. 

Cylinder renormalization has been adopted for holomorphic maps with a Siegel disk in \cite{Ya4}. We will continue with a brief summary of Yampolsky's procedure. Define ${\bf A}_\rho$ to be the Banach space of all bounded analytic functions on the open disk of radius $\rho>1$ around zero, $\field{D}_\rho \subset \field{C}$, equipped with the sup norm. ${\bf C}_\rho$ will denote the subspace of ${\bf A}_\rho$ which consists of all $f \in {\bf A}_\rho$ for which $0$ is a fixed point, $f(0)=0$, and $1$ is a critical point, $f'(1)=0$.

Given an $f \in {\bf C}_\rho$, suppose that for some $n \in \field{N}$ there exists a simple arc $l$ connecting a fixed point $a_n$ of $f^n$ with the fixed point $0$ such that $f^{-n}(l)$ is again a simple arc intersecting $l$ only at points $a_n$ and $0$. We will call the interior of the region bounded by the curves $l$ and $f^{-n}(l)$ the {\it fundamental crescent}  $C_{f^n}$  of cylinder renormalization if the iterate $f^{-n}|_{C_{f^n}}$ is univalent and the quotient of $\overline{C_{f^n} \cup f^{-n}(C_{f^n})} \setminus \{0,a_n \}$ by the iterate $f^{n}$ is conformally isomorphic to $\field{C} / \field{Z}$. 

\begin{figure}
  \begin{center}
    \setlength{\unitlength}{0.00083333in}
    \begingroup\makeatletter\ifx\SetFigFont\undefined%
    \gdef\SetFigFont#1#2#3#4#5{%
      \reset@font\fontsize{#1}{#2pt}%
      \fontfamily{#3}\fontseries{#4}\fontshape{#5}%
      \selectfont}%
    \fi\endgroup%
	{\renewcommand{\dashlinestretch}{30}
	  \begin{picture}(4930,1695)(0,-10)
	    \drawline(2022.000,445.000)(2064.009,379.174)(2116.489,321.349)
	    (2177.945,273.171)(2246.626,236.013)(2320.577,210.933)
	    (2397.693,198.644)(2475.777,199.497)(2552.606,213.467)
	    (2625.992,240.157)(2693.845,278.806)(2754.234,328.315)
	    (2805.438,387.272)(2846.000,454.000)
	    \drawline(2024.000,440.000)(2066.330,374.753)(2118.993,317.519)
	    (2180.499,269.916)(2249.110,233.289)(2322.888,208.673)
	    (2399.746,196.765)(2477.514,197.900)(2553.992,212.047)
	    (2627.020,238.806)(2694.533,277.420)(2754.623,326.798)
	    (2805.592,385.545)(2846.000,452.000)
	    \drawline(2018.000,1384.000)(2060.330,1318.753)(2112.993,1261.519)
	    (2174.499,1213.916)(2243.110,1177.289)(2316.888,1152.673)
	    (2393.746,1140.765)(2471.514,1141.900)(2547.992,1156.047)
	    (2621.020,1182.806)(2688.533,1221.420)(2748.623,1270.798)
	    (2799.592,1329.545)(2840.000,1396.000)
	    \drawline(2841.000,1396.000)(2807.023,1463.195)(2762.042,1523.580)
	    (2707.385,1575.370)(2644.668,1617.037)(2575.742,1647.349)
	    (2502.643,1665.410)(2427.531,1670.688)(2352.625,1663.026)
	    (2280.138,1642.651)(2212.209,1610.164)(2150.847,1566.526)
	    (2097.864,1513.024)(2054.825,1451.241)(2023.000,1383.000)
	    \put(4302,915){\circle{16}}
	    \put(4307,911){\circle{1176}}
	    \drawline(2022,1387)(2022,443)
	    \drawline(2845,1387)(2845,443)
	    \drawline(1285,976)(1777,976)
	    \blacken\drawline(1711.340,955.485)(1777.000,976.000)(1711.340,996.515)(1711.340,955.485)
	    \drawline(3106,964)(3599,964)
	    \blacken\drawline(3533.340,943.485)(3599.000,964.000)(3533.340,984.515)(3533.340,943.485)
	    \drawline(12,1263)(12,1262)(13,1259)
	    (14,1255)(15,1248)(17,1239)
	    (20,1227)(24,1212)(29,1195)
	    (35,1176)(42,1155)(49,1132)
	    (58,1109)(69,1085)(81,1060)
	    (94,1035)(109,1009)(126,984)
	    (146,958)(168,932)(194,905)
	    (223,878)(257,850)(294,823)
	    (336,796)(382,770)(427,747)
	    (473,727)(519,709)(565,693)
	    (610,679)(655,667)(698,656)
	    (741,647)(782,640)(824,634)
	    (864,628)(904,624)(943,620)
	    (982,617)(1019,614)(1054,612)
	    (1088,611)(1119,610)(1147,609)
	    (1172,608)(1192,608)(1210,607)
	    (1223,607)(1232,607)(1238,607)
	    (1242,607)(1243,607)
	    \drawline(4088,703)(4087,702)(4084,699)
	    (4079,694)(4072,687)(4063,677)
	    (4052,664)(4039,649)(4024,631)
	    (4009,613)(3994,592)(3979,571)
	    (3965,549)(3953,526)(3942,502)
	    (3933,477)(3926,451)(3922,424)
	    (3922,395)(3925,365)(3933,334)
	    (3946,302)(3962,274)(3981,247)
	    (4000,223)(4020,201)(4038,181)
	    (4056,164)(4071,149)(4086,135)
	    (4099,124)(4111,113)(4123,103)
	    (4135,94)(4146,85)(4159,77)
	    (4173,68)(4189,60)(4208,51)
	    (4229,43)(4254,34)(4282,27)
	    (4314,20)(4350,15)(4388,12)
	    (4428,13)(4465,17)(4501,25)
	    (4535,34)(4566,44)(4594,55)
	    (4620,66)(4643,77)(4663,88)
	    (4681,98)(4698,108)(4713,117)
	    (4726,127)(4740,136)(4753,146)
	    (4765,156)(4778,167)(4792,180)
	    (4806,193)(4820,209)(4836,226)
	    (4851,245)(4867,267)(4882,291)
	    (4896,317)(4907,344)(4915,373)
	    (4918,404)(4917,435)(4911,464)
	    (4901,492)(4888,517)(4872,542)
	    (4854,565)(4833,587)(4812,608)
	    (4788,628)(4764,647)(4739,665)
	    (4715,683)(4690,699)(4667,714)
	    (4646,727)(4627,738)(4611,748)
	    (4598,755)(4578,767)
	    \blacken\drawline(4644.858,750.810)(4578.000,767.000)(4623.748,715.627)(4644.858,750.810)
	    \drawline(12,1263)(13,1263)(17,1264)
	    (23,1265)(32,1267)(44,1270)
	    (60,1273)(80,1276)(102,1279)
	    (127,1283)(155,1286)(184,1289)
	    (215,1291)(247,1291)(281,1291)
	    (316,1289)(352,1285)(391,1280)
	    (431,1271)(474,1260)(520,1246)
	    (568,1228)(618,1207)(669,1181)
	    (713,1156)(755,1130)(795,1103)
	    (832,1075)(868,1047)(901,1019)
	    (932,991)(961,963)(989,935)
	    (1015,907)(1040,879)(1063,851)
	    (1086,824)(1107,797)(1128,770)
	    (1147,745)(1164,720)(1181,698)
	    (1195,678)(1207,660)(1218,644)
	    (1227,632)(1233,622)(1238,615)
	    (1241,611)(1242,608)(1243,607)
	    \drawline(437,1014)(435,1014)(430,1014)
	    (421,1014)(408,1015)(391,1015)
	    (371,1015)(348,1014)(324,1013)
	    (299,1011)(273,1007)(248,1003)
	    (222,997)(197,989)(172,979)
	    (147,965)(124,949)(102,929)
	    (86,909)(72,889)(62,869)
	    (53,851)(47,835)(42,820)
	    (38,807)(35,795)(33,785)
	    (31,774)(30,764)(29,753)
	    (29,742)(30,730)(33,716)
	    (36,700)(42,684)(51,666)
	    (63,649)(78,632)(97,618)
	    (118,607)(139,598)(159,591)
	    (179,586)(196,583)(213,580)
	    (228,578)(242,577)(256,576)
	    (269,576)(284,576)(299,577)
	    (316,578)(335,580)(355,584)
	    (378,590)(402,597)(426,608)
	    (449,621)(472,639)(490,660)
	    (504,682)(514,704)(522,727)
	    (527,751)(529,775)(530,799)
	    (530,822)(529,845)(528,866)
	    (526,885)(524,901)(520,929)
	    \blacken\drawline(549.595,866.901)(520.000,929.000)(508.977,861.099)(549.595,866.901)
	    \drawline(2520,694)(2522,692)(2525,689)
	    (2531,682)(2541,673)(2554,661)
	    (2569,645)(2588,628)(2608,609)
	    (2630,590)(2653,570)(2678,550)
	    (2703,531)(2729,513)(2756,496)
	    (2785,480)(2815,466)(2847,454)
	    (2880,445)(2914,439)(2946,437)
	    (2976,438)(3002,442)(3025,446)
	    (3043,452)(3059,457)(3072,463)
	    (3083,469)(3092,474)(3100,481)
	    (3108,487)(3116,494)(3124,502)
	    (3133,512)(3143,523)(3152,536)
	    (3163,552)(3172,571)(3180,592)
	    (3185,615)(3186,639)(3183,663)
	    (3178,685)(3172,705)(3165,723)
	    (3158,738)(3152,751)(3145,763)
	    (3139,773)(3133,783)(3126,793)
	    (3118,803)(3109,814)(3098,825)
	    (3084,839)(3068,853)(3048,869)
	    (3025,886)(2999,903)(2969,918)
	    (2934,931)(2899,941)(2866,946)
	    (2833,949)(2802,949)(2772,946)
	    (2743,942)(2715,937)(2688,931)
	    (2663,924)(2639,917)(2619,910)
	    (2601,905)(2572,894)
	    \blacken\drawline(2626.116,936.468)(2572.000,894.000)(2640.668,898.105)(2626.116,936.468)
	    \put(1491,1056){\makebox(0,0)[lb]{\smash{{{\SetFigFont{7}{8.4}{\rmdefault}{\mddefault}{\updefault}$\Phi$}}}}}
	    \put(1283,621){\makebox(0,0)[lb]{\smash{{{\SetFigFont{7}{8.4}{\rmdefault}{\mddefault}{\updefault}$a_n$}}}}}
	    \put(9,1311){\makebox(0,0)[lb]{\smash{{{\SetFigFont{7}{8.4}{\rmdefault}{\mddefault}{\updefault}$0$}}}}}
	    \put(681,539){\makebox(0,0)[lb]{\smash{{{\SetFigFont{7}{8.4}{\rmdefault}{\mddefault}{\updefault}$l$}}}}}
	    \put(796,1160){\makebox(0,0)[lb]{\smash{{{\SetFigFont{7}{8.4}{\rmdefault}{\mddefault}{\updefault}$f^{-n}(l)$}}}}}
	    \put(85,454){\makebox(0,0)[lb]{\smash{{{\SetFigFont{7}{8.4}{\rmdefault}{\mddefault}{\updefault}$R_{f^n}$}}}}}
	    \put(4283,969){\makebox(0,0)[lb]{\smash{{{\SetFigFont{7}{8.4}{\rmdefault}{\mddefault}{\updefault}$0$}}}}}
	    \put(3306,1013){\makebox(0,0)[lb]{\smash{{{\SetFigFont{7}{8.4}{\rmdefault}{\mddefault}{\updefault}$e$}}}}}
	    \put(2898,1234){\makebox(0,0)[lb]{\smash{{{\SetFigFont{7}{8.4}{\rmdefault}{\mddefault}{\updefault}$\field{C} / \field{Z}$ }}}}}
	    \put(2800,297){\makebox(0,0)[lb]{\smash{{{\SetFigFont{7}{8.4}{\rmdefault}{\mddefault}{\updefault}$\Phi \circ R_{f^n} \circ \Phi^{-1}$}}}}}
	    \put(4220,151){\makebox(0,0)[lb]{\smash{{{\SetFigFont{7}{8.4}{\rmdefault}{\mddefault}{\updefault}${}_n R_{cyl}(f)$}}}}}
	  \end{picture}
	}
	\label{cyl_renorm} \caption{Schematics of the cylinder renormalization.}
  \end{center}
\end{figure}
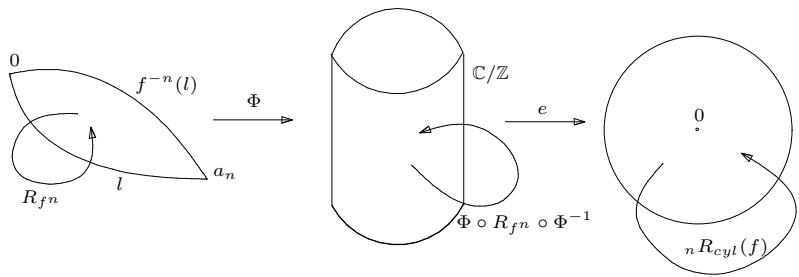

Furthermore, let us assume that $R_{f^n}$, the first return map of $C_{f^n}$, has a critical point $c \in C_{f^n}$.  The conformal isomorphism of $C_{f^n}$ and $\field{C} / \field{Z}$, normalized so that it maps $c$ to $1$, will be denoted by $\Phi$. The map ${}_n R_{cyl}(f) = e \circ \Phi \circ R_{f^n} \circ \Phi^{-1} \circ e^{-1}$ (see Fig. $1$), where $e(z)=e^{-2 \pi i z}$, will be called a {\it cylinder renormalization} of $f$ with period $n$ if ${}_n R_{cyl}(f) \in {\bf  C}_{\rho'}$ for some $\rho'>1$. Naturally, this definition should be independent of the particular choice of the fundamental crescent. In fact, any other fundamental crescent $C'_{f^{n}}$ with the same endpoints as  $C_{f^{n}}$, and such that  $C_{f^{n}} \cup C'_{f^{n}}$ is a topological disk, produces the same ${}_n R_{cyl}(f)$ (cf \cite{Ya}).

It is clear from this definition that the knowledge of the conformal isomorphism of $C_{f^n}$ and $\field{C}/\field{Z}$ is of a key importance for building a cylinder renormalization of a map. In the remaining part of the paper we will describe how to construct such an isomorphism for an analytic map of degree $d \ge 2$  with an irrationally indifferent fixed point at $0$:
\begin{equation}
\nonumber f(0)=0 \ {\rm and}  \ f'(0)=e^{2 \pi i \theta}.
\end{equation} 

We will further assume that $\theta$ is of {\it bounded type}, that is it admits a  continued fraction expansion
\begin{equation} \label{cont_frac}
\theta={1 \over r_0+{1 \over r_1+{1 \over r_2+\ldots}}  },
\end{equation}
(which we will abbreviate as $\theta=[r_0,r_1,r_2,\ldots]$) with 
\begin{equation}
\nonumber \sup r_i < \infty.
\end{equation}

 The rational number corresponding to the $n$-th truncation of this expansion is commonly denoted by $q_n/p_n$: 
\begin{equation}
\nonumber {q_n \over p_n}=[r_0,r_1,\ldots,r_n].
\end{equation}

By the classical result of \citeasnoun{Si}, for such an $f$ there exists a maximal open connected neighborhood $U$ of $0$, called the {\it Siegel disk}, on which the action of $f$ is conformally conjugate to the rigid rotation $z \mapsto e^{2 \pi i \theta} z$ of the open unit disk. 

We will now focus on the construction of a cylinder renormalization for analytic maps with a Siegel disk. X. Buff (unpublished) has suggested a construction of a fundamental domain for a quadratic polynomial as the interior of a region bounded by internal and external rays, logarithmic spirals and  equipotentials. However, the exponential density of the Julia set makes a computer implementation of some of these curves practically difficult. 

Here, we will make a different choice of a fundamental crescent for functions with a Siegel disk. The curve $l_n$ will consist of two parabolas: one passing through points $0$ and $f^{q_{n+2}+q_n}(1)$, the other -- through point $f^{q_{n+2}+q_n}(1)$ and a repelling fixed point $a_{q_n}$. Such two parabolas are uniquely defined after one specifies their common tangent line at point $f^{q_{n+2}+q_n}(1)$. The slope of this line can be chosen in a convenient way. Of course our, choice of the boundary $l_n$ is rather arbitrary, but has the virtue of having a simple analytic form. This simple choice of boundaries can be used for a general analytic function in ${\bf C}_\rho$, and not only for quadratic polynomials.

We shall take the interior of the region bounded by the curves $l_n$ and $f^{-q_n}(l_n)$ as a candidate for a fundamental crescent $C_{f^{q_n}}$ (here we chose the inverse branch of $f^{-q_n}$ that maps $0$ into itself). One has to verify that the simple curve $f^{-q_n}(l_n)$ intersects  $l_n$ only at the endpoints. We will momentarily defer the explanation of this verification.

As we have already explained, the key step in the cylinder renormalization of a function $f$ is a construction of a properly normalized conformal isomorphism $\Phi$.
It should be noted however that this uniformizing change of coordinates can not be written down explicitly, and, as we will show, requires a solution of a Beltrami equation which, in turn, can not be found ``by hand''. One, therefore, may try to find $\Phi$ numerically and provide bounds in this numerical solution: Below, we will demonstrate how this can be done for the golden-mean quadratic polynomial $f_{\theta^*}$, and $n=3$.

\medskip   \section{ Uniformization of a fundamental crescent of cylinder renormalization}\label{modelproblem}
\setcounter{equation}{0}

\begin{figure}
  \begin{center}
    \resizebox{120mm}{!}{\includegraphics{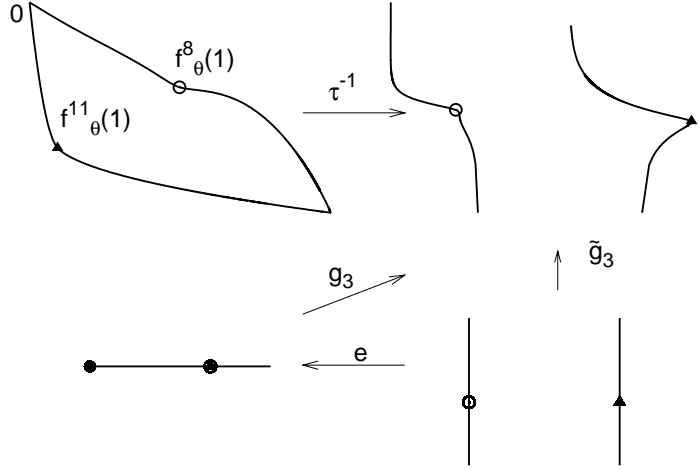}}
    \label{schematics}    \caption{Mapping $g^{-1}_1 \circ \tau^{-1} $ of the interior of the fundamental crescent $C_{f_{\theta^*}}$ onto the punctured plane with a slit $\field{C}^* \setminus \field{R}^+$.}
  \end{center}
\end{figure}

In this Section we will outline the steps in the construction of the conformal isomorphism $\Phi$  (see Fig. ~\ref{schematics} for the schematics of our procedure). Most of the details of this construction are independent of the choice of $f$, however at the end of this Section we will specify our fundamental crescent for the golden-mean quadratic polynomial and $n=3$.

Introduce  for all $z \in C_{f^{q_n}}$ a new coordinate $\xi$ through
\begin{equation}
  z=\tau(\xi)={a_{q_n} \over (1-e^{i a  \xi+b })}, 
\end{equation}
where normalizing constants $a$ and $b$ are chosen so that $|\tau^{-1}(f^{q_{n+2}+q_n}(1))|=1$ and $\tau^{-1}(f^{q_{n+2}}(1))=0$. The choice of of this coordinate is motivated by the fact that $\tau^{-1}$ maps the interior of the fundamental crescent conformally onto the interior of an infinite vertical closed strip $\mathcal{S}$, whose width is comparable to one. This coordinate has been suggested in \cite{Shish} in a similar context. 

Next, following \cite{Shish},  define a function $\tilde{g}_n$ from $\mathcal{U}\overset{{\rm def}}=\left\{u+i v \in \field{C}\right.: 0 \le \Re {w} \le 1\left. \!\!\!\!\!\!\phantom{C}\right\}$ to $\mathcal{S}$ by setting
\begin{align}
  \nonumber   \tilde{g}_n(u+i v)&= (1-u) \tau^{-1} (f^{-q_n}(\gamma_n(v)))+u \tau^{-1} (\gamma_n(v))  
\end{align}
(below, we will specify the choice of parametrization $\gamma_n : \field{R}  \mapsto  l_n$). Let  $\sigma_0$  be the standard conformal structure on $\field{C}$, and let $\sigma=\tilde{g}_n^* \sigma_0$ be its pull-back on $\mathcal{U}$. Extend this conformal structure to $\field{C}$ through $\sigma=(T^k)^* \sigma$ on  $T^{-k}(\mathcal{U})$,  $T(w)=w+1$.   By the Measurable Riemann Mapping theorem (see Section \ref{Beltrami}) there exists a unique quasiconformal mapping $\tilde{g}: \field{C} \mapsto  \field{C}$ such that $\tilde{g}^* \sigma_0 = \sigma$, normalized so that $\tilde{g}(0)=0$ and $\tilde{g}(1)=1$. Notice that $\tilde{g} \circ T \circ \tilde{g}^{-1}$  preserves the standard conformal structure:
\begin{equation}
\nonumber \left(\tilde{g} \circ T \circ \tilde{g}^{-1}\right)^* \sigma_0 = (\tilde{g}^{-1})^* 
\circ T^* \circ \tilde{g}^* \sigma_0 =(\tilde{g}^{-1})^* \circ T^* \sigma=(\tilde{g}^{-1})^* \sigma =(\tilde{g}^*)^{-1} \sigma=\sigma_0,
\end{equation}
therefore it is a conformal automorphism of $\field{C}$, i.e. an affine map, and in fact, a translation, since $T$ does not have fixed points in $\field{C}$. Moreover $\tilde{g} \circ T \circ \tilde{g}^{-1}(0)=1$, Therefore  $\tilde{g} \circ T \circ \tilde{g}^{-1}=T$.

By the definition of $\tilde{g}_n$, $\tilde{g}_n^{-1} \circ \tau^{-1} \circ f^{q_n} \circ \tau=T \circ \tilde{g_n}^{-1}$ on the image of $l_n$ by $\tau^{-1}$. Therefore, the map $\phi=\tilde{g} \circ \tilde{g}_n^{-1}$, defined on $\mathcal{S}$, can be continuously extended to all of $\field{C}$. Since $\phi$ is conformal on the interior of $\mathcal{S}$, $\oint_\gamma \phi(w) d w=0$ over any closed contour $\gamma$ in $\field{C}$, and by Morera's theorem $\phi$ is analytic on all of $\field{C}$. Clearly, $\Phi=\phi \circ \tau^{-1}$ is the desired conformal isomorphism of $C_{f^{q_n}}$ and $\field{C} / \field{Z}$.  

Next, define $g=e \circ \tilde{g} \circ e^{-1}$. Since 
\begin{equation}
\nonumber {g_{\bar{z}}(e(w)) \over g_z(e(w))}={e(w) \over \overline{e(w)}} {\tilde{g}_{\bar{w}}(w) \over \tilde{g}_w(w)},  
\end{equation}
the $1$-periodic function $\tilde{g}$ is a solution of the {\it Beltrami equation} 
\begin{equation}
\nonumber \tilde{g}_{\bar{w}}=\tilde{\mu}\tilde{g}_{w}, \quad \tilde{\mu}=(\tilde{g}_n)_{\bar{w}}/(\tilde{g}_n)_{w}
\end{equation}
whenever $g$ is a solution of  
\begin{equation}
g_{\bar{z}}=\mu g_{z}, \quad \mu(z)=(z/\bar{z}) \tilde{\mu}(e^{-1}(z)).
\end{equation}

Thus, we have reduced the problem of finding $e \circ \Phi=g \circ e \circ \tilde{g}^{-1}_n \circ \tau^{-1}$ to that of finding the properly normalized solution of the Beltrami equation
\begin{equation}\label{beltt}
g_{\bar{z}}=\mu g_z, \quad \mu(z)={z \over \bar{z}} {(\tilde{g}_n)_{\bar{w}} (e^{-1}(z)) \over (\tilde{g}_n)_{w} (e^{-1}(z)) } 
\end{equation}
on the punctured plane $\field{C}^*$.  The issue of existence of a solution of $(\ref{beltt})$ is addressed by the famous Measurable Riemann Mapping theorem, also known as Ahlfors-Bers-Boyarskii theorem (see \cite{AB} and \cite{Boyarskii}):
\begin{thm} \label{ABB}(Ahlfors and Bers, Boyarskii).
Let $\mu \in L_\infty(\hat{\field{C}})$ satisfy  $\left\| \mu \right\|_\infty \le K < 1$. Then there exists a unique map $g^\mu: \hat{\field{C}} \rightarrow \hat{\field{C}}$ such that $g^\mu_{\bar{z}} (z)/ g^\mu_z(z)=\mu(z)$, $g^\mu$ fixes $0$, $1$ and $\infty$, and $g^\mu$ is a $(1+K)/(1-K)$-quasiconformal map. 
\end{thm}


Recall, that there are two equivalent definitions of a quasiconformal map, the ``geometric'' and the ``analytic'' one. Here we will give only the ``analytic'' definition which is shorter: an interested reader is referred to \cite{Ahlfors}, \cite{Bers}, \cite{Markovic} and \cite{Boyar_Iwan} for a fuller account of the basics of the theory.

Given an open set  $\Omega \subset \field{C}$, a map $f: \Omega \rightarrow f(\Omega) \subset \field{C}$ is said to be absolutely continuous on lines (ACL) in a rectangle $R \subset \Omega$ with sides parallel to the $x$ and $y$ axes if $f$ is absolutely continuous on almost every horizontal and vertical line in $R$. The map $f$ is ACL on $\Omega$, if it is ACL on every rectangle in $\Omega$. Partial derivatives $f_z=(f_x-i f_y)/2$ and $f_{\bar{z}}=(f_x+i f_y)/2$ of such map exist a.e.  in $\Omega$.
\begin{defn}
A homeomorphism  $f: \Omega \rightarrow f(\Omega)$  is K-quasiconformal if and only if the following holds:
\begin{itemize}
\item[i)] $f$ is ACL on every rectangle in $\Omega$,
\item[ii)] $|f_{\bar{z}}| \le {K-1 \over K+1} |f_z|$ a.e. in $\Omega$.
\end{itemize}
\end{defn}

The complex dilatation of a quasiconformal homeomorphism $f$ is $\mu_f(z)=f_{\bar{z}}(z)/f_z(z)$.  A mapping $f$ with a prescribed complex dilatation $\mu$ solves the following Beltrami equation: 
\begin{equation}\label{belt_eqn}
f_{\bar{z}}(z)=\mu(z)  f_z(z),
\end{equation}
and its existence is the subject of Theorem \ref{ABB}.

\medskip   \section{Statement of main results}\label{Results}
\setcounter{equation}{0} 

Most of the present work deals with a constructive  implementation of the original  proof of the Measurable Riemann Mapping Theorem due to L. Ahlfors and L. Bers. One of our main objectives will be to present rigorous bounds on how far a numerical solution of $(\ref{beltt})$ lies from the actual one. We will give a detailed derivation of the necessary  bounds in Sections \ref{CZconstant}--\ref{Bounds}. At the moment, however, we will proceed to summarize our results.

\begin{thm}\label{Main_Theorem}
Let $\mu \in L_\infty(\hat{\field{C}})$ and an integer $p>2$  be such that $\left\| \mu \right\|_\infty\le K < 1$ and $K C_p<1$, where 
\begin{equation}
\nonumber C_p= \cot^2(\pi/2p).
\end{equation}
 Assume that $\mu=\nu+\eta+\gamma$, where $\nu$ and $\eta$ are compactly supported in $\field{D}_R$, and  $\gamma(z)$ is supported in $\hat{\field{C}} \setminus \field{D}_R$.  Furthermore, let $\eta$ be in $L_p(\field{D}_R)$ and $\left\| \eta \right\|_p < \delta$ for some sufficiently small $\delta$. Also, let $h^* \in L_p(\field{C})$ and $\epsilon$   be such that $B_p(h^*,\epsilon)$, the ball of radius $\epsilon$ around $h^*$ in  $L_p(\field{C})$, contains  $B_p(T_\nu[h^*],C_p \epsilon')$, with 
\begin{equation}
\nonumber \epsilon'= \delta \  \esssup_{\field{D}_R}{ |h^*  +1|}+K  \epsilon.
\end{equation}
Then the solution $g^\mu$ of the Beltrami equation $g^\mu_{\bar{z}}=\mu g^\mu_z$ admits the following bound:
\begin{equation}\label{mainbound}
\left|g^\mu(z)-g^\nu_*(z)\right|\le   A \epsilon' |z|^{1-2/p}  + {C \left[ \left|g^\nu_*(z)\right|+A \epsilon' |z|^{1-2/p} \right]^{1+2/p} \over R^{4/p} -C \left[ \left|g^\nu_*(z)\right|+A \epsilon' |z|^{1-2/p} \right]^{2/p}},
\end{equation}
where $g_*^\nu(z)=P_\nu[h^*](z)+z$ and 
\begin{align}
\nonumber A&={1 \over \pi^{1/p}} \left[ 4 { 3 ^{2 q-2} \over 2^q (2-q)} +{25 \over 36} 3^{2 q}-23^{2q-2}+ {\left( {4 \over 9}\right)^{1-q} \over q-1 } \right]^{1 / q},\\
\nonumber C&=\pi^{1/p} A K { R^{4/p} ( 2-K C_p ) \over  r^{2/p}( 1-K C_p)},\quad r=\inf_{\left\{ z \in \field{C}: |z|=R \right\}}\left\{\left| g^{\nu+\eta}(z)  \right| \right\}.
\end{align}
\end{thm}

We have used the bounds from Theorem $\ref{Main_Theorem}$  to uniformize a fundamental domain for the cylinder renormalization of period $3$ of the golden mean quadratic polynomial. Before we state our findings, we would first like to specify our choice of a fundamental crescent $C_{f^3}$.

It is convenient to parametrize the boundaries  $l_n$ and $f^{-q_n}(l_n)$ of $C_{f^{q_n}}$ by the radial coordinate in $\field{C}$. Of course, this parametrization  is far from unique. Here, we will mention our choice for $f_{\theta^*}(z)=e^{2 \pi i \theta^*} z(1-0.5 z)$ and $n=3$, motivated by our numerical experiments: 
\begin{equation}\label{l_n}
\lambda_3(r) = \left\{ (x(r),A x(r)^2+B x(r)), \quad r \le \tilde{r},  ~ \atop  (C y(r)^2+D y(r)+E,y(r)) , \quad r>\tilde{r},  \right. 
\end{equation}
where 
\begin{align}
  \nonumber x(r)&={\Re{f^{11}_{\theta^*}(1)} \over |f^{11}_{\theta^*}(1)  |} T(r), \\  
  \nonumber y(r)&=\Im{f^{11}_{\theta^*}(1)} {|a_3-f^{11}_{\theta^*}(1)|+|f^{11}_{\theta^*}(1)|-T(r) \over |a_3-f^{11}_{\theta^*}(1)|}+\Im{a_3} {T(r)-|f^{11}_{\theta^*}(1)|) \over |a_3-f^{11}_{\theta^*}(1)|},\\
  \nonumber T(r)&= {|a_3-f^{11}_{\theta^*}(1)|+|f^{11}_{\theta^*}(1)| \over \sqrt{r} +1} \sqrt{r},
\end{align}
and $\tilde{r}$ is defined through $T(\tilde{r})=|f^{11}_{\theta^*}(1)|$. Constants $A$, $B$, $C$, $D$ and $E$ are fixed by the conditions $0, f^{11}_{\theta^*}(1), a_3 \in  l_3$, together with the requirement that the slope of the  common tangent line to both parabolas at point  $f^{11}_{\theta^*}(1)$ is equal to $-1.7$.

Define the following function on $\field{C}^*$:
\begin{equation}
 \nonumber   g_n(r,\phi)= \left(\eta(-\phi)+{\phi \over 2 \pi}\right) \tau^{-1}(f^{-q_n}(\lambda_n (r)))+\left( 1-\eta(-\phi)- {\phi \over 2 \pi} \right)\tau^{-1}(\lambda_n (r)),
\end{equation}
where $-\pi < \phi \le \pi$, and $\eta$ is the Heaviside step function (we have adopted the convention $\eta(0)=1$). It is clear that $\tilde{g}_n \circ e^{-1}(r e^{ i \phi})=g_n(r,\phi)$. Therefore if $\mu$ is as in $(\ref{beltt})$, then
\begin{equation}\label{belt2}
\mu(r e^{i\phi})=e^{2 i \phi} {r \partial_r g_n(r,\phi)+i \partial_{\phi} g_n(r,\phi) \over r \partial_r g_n(r,\phi)-i \partial_{\phi} g_n(r,\phi) } 
\end{equation}
on $\field{C}^*$.

One also needs to verify the empty intersection of the boundaries of the fundamental domain. To this end, choose a sufficiently large $R$ and a sufficiently small $\varrho$,  and consider the arc $\lambda_n([\varrho,R])$. Cover this arc by a finite collection $\left\{B_i\right\}^{N}_{i=1}$ of open balls in $\field{C}$ (see Appendix for a discussion of evaluation of analytic functions on open disks, referred to as ``standard sets in $\field{C}$) and verify that $f^{-q_n}\left( B_k \right) \cup \left\{B_i\right\}^{N}_{i=1}=\emptyset$, $k=1 \ldots N$. Verification of the empty intersection of the boundary curves for small, $r<\varrho$, and large, $r>R$, values of the parameter is not needed: Rather than using parabolas, we will complete the boundaries of the fundamental crescent for these parameter values with pieces of rays and logarithmic spirals together with their preimages (see Section \ref{Completion} for details). We will refer to this construction simply as a ``completion'' of the fundamental crescent. 
\begin{figure}
  \begin{center}
\begin{tabular}{c c}   
 \resizebox{60mm}{!}{\includegraphics{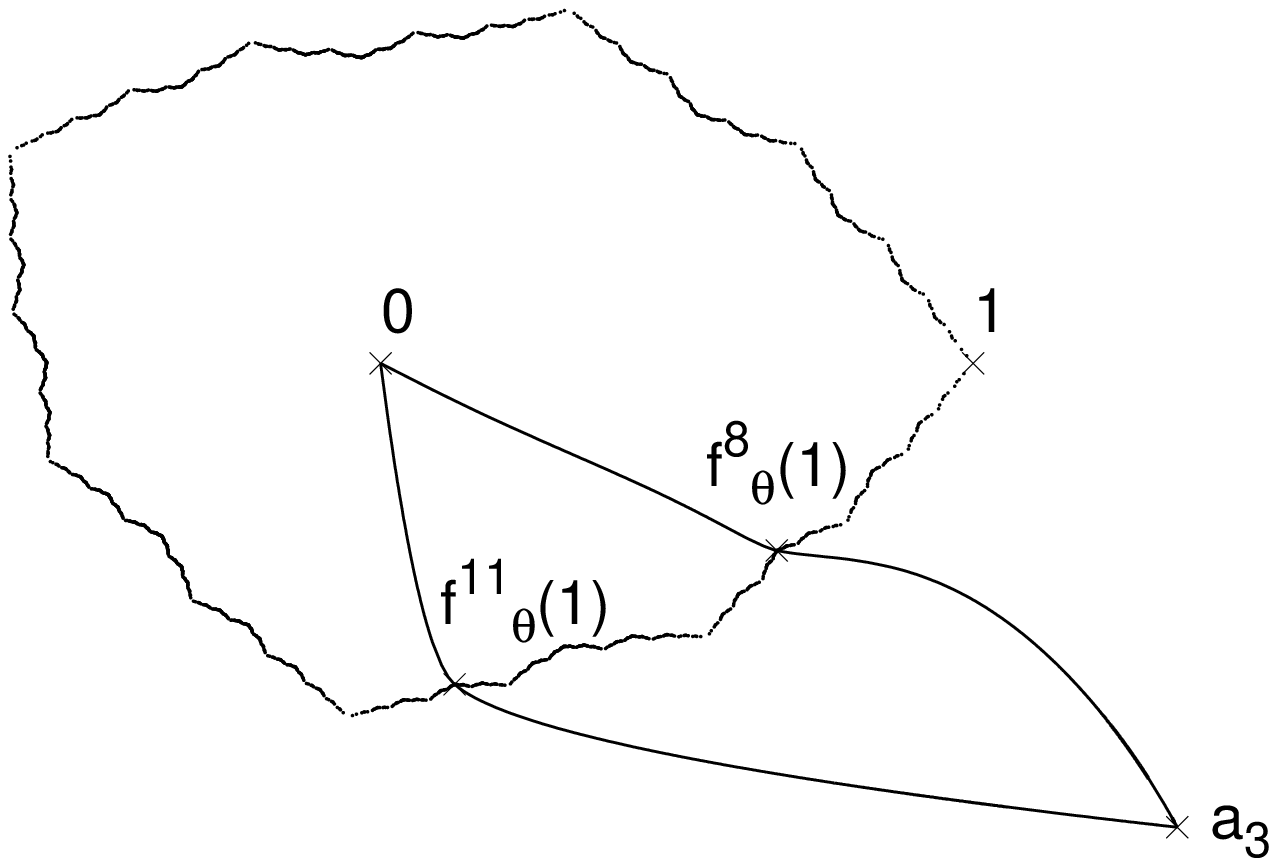}} & \resizebox{50mm}{!}{\includegraphics{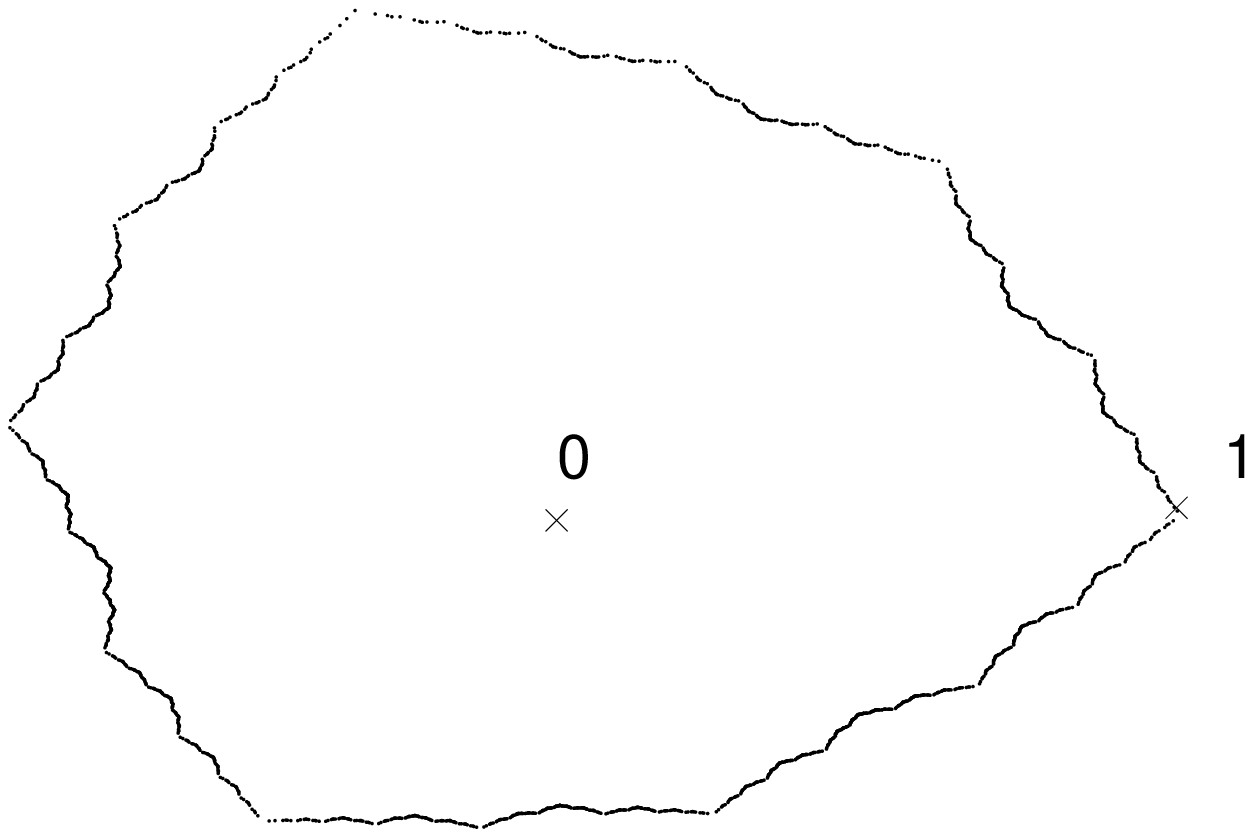}}\\
(a) & (b) 
\end{tabular}
    \caption{The fundamental domain for the cylinder renormalization of the quadratic polynomial $f_{\theta^*}(z)=e^{2 \pi i {\theta^*}} z(1-0.5 z)$, ${\theta^*}=(\sqrt{5}-1)/2$, of period 3 (a) and the Siegel disk of the cylinder renormalization of $f_{\theta^*}$ of period 3 (b).}
    \label{test}
  \end{center}
\end{figure}
We have performed the  verification of the empty intersection for the golden mean quadratic polynomial $f_{\theta^*}(z)=e^{2 \pi i {\theta^*}} z(1-0.5 z)$, and $n=3$.  Although most of the results of this paper are valid for any analytic $f \in {\bf C}_\rho$, our Main Example, that serves purely demonstrative purposes, will be stated only for these particular data. 
  
\begin{mainexample}\label{mainex}
Let $f_{\theta^*}(z)=e^{2 \pi i {\theta^*}} z(1-0.5 z)$, where ${\theta^*}=(\sqrt{5}-1)/2$, and let $\lambda_3$ be as in 
$(\ref{l_n})$. Then, there exists a completion of the fundamental crescent $C_{f^3_{\theta^*}}$ bounded by the arcs 
$\lambda_3([\varrho,R])$ and $f^{-3}_{\theta^*}( \lambda_3 ([\varrho,R]) )$ with $\varrho=3.1 \times 10^{-3}$ and 
$R=9.8 \times 10^3$, and a finite Fourier series $g_*$ whose dilatation is compactly supported in $\field{D}_R$ such 
that the solution $g^\mu$ of the Beltrami equation  $(\ref{beltt})$ on $\hat{\field{C}}$ satisfies
\begin{equation}
\nonumber \left|g^\mu(z)-g_*(z)\right|\le   A \epsilon' |z|^{1-2/p}  + {C \left[ \left|g_*(z)\right|+A \epsilon' |z|^{1-2/p} \right]^{1+2/p} \over R^{4/p} -C \left[ \left|g_*(z)\right|+A \epsilon' |z|^{1-2/p} \right]^{2/p}},
\end{equation}
for $p=2.1$ and $A <7.11$, $\epsilon' < 2.66$ and $C<6.14$.
\end{mainexample}

All bounds on the constants reported in the preceding Example have been found with the help of a computer (see the following Sections for a detailed discussion of the ingredients of the computer-assisted proof). To this end, we have written a set of routines in the programming language Ada 95 (cf \cite{ADA} for the language standard). We have parallelized our programs and compiled them with the public version 3.15p of the GNAT compiler \cite{GNAT}. The programs \citeyear{programs} have been run on the computational cluster of 64 2.2 GHz AMD Opteron processors located at the University of Texas at Austin.

We believe that out constructive bounds on the solution of the Beltrami equation arising in cylinder renormalization are a major step towards rigorous computer-assisted study of this operator. Apart from that, a constructive proof of the Measurable Riemann Mapping Theorem is interesting in its own right, and should find other applications in one-dimensional complex dynamics and geometry.

\medskip   \section{A constructive proof of the Measurable Riemann Mapping Theorem}\label{Beltrami}
\setcounter{equation}{0} 

 The classical proof of Theorem $\ref{ABB}$ uses two integral operators in an essential way; the first operator being the Hilbert transform
\begin{equation}
\nonumber T [h](z)={i \over 2 \pi} \lim_{\epsilon \to 0} \int \int_{\field{C} \setminus B(z,\epsilon)} {h(\xi) \over (\xi -z)^2} \ d \bar\xi \wedge d \xi, 
\end{equation}
the second --- the Cauchy transform
\begin{equation}
\nonumber P [h] (z)={i \over 2 \pi} \int \int_{\field{C}} {h(\xi) \over (\xi -z)} \ d \bar\xi \wedge d \xi.
\end{equation}

 Some of the properties of the Hilbert transform are described in the celebrated Calderon-Zygmund lemma (see \cite{CalZyg}):
\begin{lemma}(Calderon and Zygmund).
The Hilbert transform is a well-defined  bounded operator on $L_p(\field{C})$ for all $2<p<\infty$: for every such $p$ there exists a constant $C_p$ such that $\left\| T [h] \right\|_p \le C_p \left\| h \right\|_p$ for any $h \in L_p(\field{C})$, and $C_p \rightarrow 1$ as $p \rightarrow 2$. 

\end{lemma}

For compactly supported differentials one can state an extended version of the Ahlfors-Bers-Boyarskii theorem:
\begin{thm} \label{compactABB}
Let $\mu \in L_\infty(\field{C})$ be compactly supported and satisfy $\left\| \mu \right\|_\infty \le K < 1$. 
Then for every $p>2$ such that $C_p K<1$, the operator $h \rightarrow T [ \mu (h+1)] $ is a contraction on $L_p(\field{C})$ with a unique fixed point $h^*$. Moreover, the solution of the Beltrami equation, $g^\mu$, is  given by 
\begin{equation}\label{solution}
 g^\mu=P [\mu (h^*+1)] +id,
\end{equation}
and is such that $g(0)=0$, $g$ is continuous, has distributional derivatives, and $g_z-1 \in L_p(\field{C})$.

\end{thm}

In the rest of the paper the operators $T [ \mu (\cdot+1)]$ and $P [ \mu (\cdot+1)] $ will be abbreviated as $T_\mu$ and $P_\mu$, respectively.

The above version of the  Ahlfors-Bers-Boyarskii theorem for compactly supported differentials provides the building blocks for the standard proof of Theorem $\ref{ABB}$ (see \cite{Ahlfors}). Here we will present a brief sketch of this proof. 

Let $\mu=\nu+\gamma$ where $\gamma(z)=0$ whenever $z \in \field{D}_R$, the disk of the radius $R$ around $0$, and $\nu(z)=0$ whenever $z \in \hat{\field{C}} \setminus \field{D}_R$.  Let $g^{\nu}$ be the solution of the Beltrami equation $g^{\nu}_{\bar{z}}=\nu g^{\nu}_z$. First, define $\beta$ by setting
\begin{equation}\label{beta}
\beta \circ g^{\nu}= \gamma {g^{\nu}_z \over \bar{g}^{\nu}_{\bar{z}} }.  
\end{equation}
Next, one obtains the solution $g^{\beta}$ of $g^{\beta}_{\bar{z}}=\beta g^{\beta}_z$: Define 
\begin{equation}\label{betatilde}
\tilde{\beta}(z)=\beta(R^2/z) z^2/\bar{z}^2,
\end{equation}
which is compactly supported in a neighborhood of zero, and let $g^{\tilde{\beta}}$ solve the Beltrami equation for $\tilde{\beta}$. Then, $g^{\beta}(z)=R^2/g^{\tilde{\beta}}(R^2/z)$ solves $g^\beta_{\bar{z}}=\beta g^\beta_z$. Finally, $g^{\mu}=g^\beta \circ g^{\nu}$ solves $g^\mu_{\bar{z}}=\mu g^\mu_z$.

The Beltrami equation arises naturally in problems of uniformization of domains in $\field{C}$. Many uniformization problems require rigorous estimates on the solution of the Beltrami equation. In some cases one can obtain rough, but sufficient estimates ``by hand'' (cf \cite{Shish} for an example similar to one considered here); often, however, one requires much tighter bounds on the solution.  Although there are efficient numeric algorithms for a Beltrami equation with a compactly supported $\mu$ (cf \cite{Daripa1}, \cite{Daripa2}, \cite{Daripa3} and \cite{GaiKhmel}), rigorous estimates on how far such numeric solutions are from a true solution have not been yet provided in literature. Our main objective will be to describe how one can obtain such estimates on an approximate numerical solution. In doing this, we will consider a Beltrami equation for the differential $(\ref{belt2})$ that arises in the abovementioned problem of cylinder renormalization.

\medskip\section{The Calderon-Zygmund constant}\label{CZconstant}
\setcounter{equation}{0}

In our following exposition we will require an estimate on the Calderon-Zygmund constant.
\begin{lemma}
The Calderon-Zygmund constant $C_p$ satisfies 
\begin{equation}
\nonumber C_p \le  \cot^2(\pi/2p).
\end{equation}
\end{lemma}
\begin{proof}
For $f(x) \in C^1_0(\field{R})$ define the one-dimensional Hilbert transform
\begin{equation}
\nonumber H[f](x)=P.V. {1 \over \pi} \int^{\infty}_{-\infty}{f(y) \over y-x} d y.  
\end{equation}

Application of the Calderon-Zygmund lemma to this singular transform yields 
\begin{equation}\label{1DCalZyg}
\left\| H[f] \right\|_p<A_p \left\| f\right\|_p.
\end{equation}

As it has been shown in \cite{Grafakos}, the best bound on $A_p$ for $2 \le p < \infty$ is  $\cot(\pi/2p)$.
Next, for $g \in C_0^2(\field{C})$ define
\begin{equation}
\nonumber T^*[g](z)=-{ i \over 4 \pi} \int \int_{\field{C}} {g(\xi+z) \over \xi |\xi| } d \bar{\xi} \wedge d \xi.
\end{equation}
By a classical result (cf \cite{Ahlfors}) $T[g]=-T^*[T^*[g]]$ for $g \in C^2_0(\field{C})$, and this equality can be extended to $L_p(\field{C})$. Moreover,  $T^*[g] \in L_p(\field{C})$ whenever $g \in L_p(\field{C})$,  and 
\begin{equation}
\nonumber T^*[g](z) ={1 \over 2} \int_0^\pi \left( {1 \over \pi} \int_0^{\infty} {g(z+r e^{i \phi})-g(z-r e^{i \phi})  \over r } d r \right) e^{-i \phi} d \phi. 
\end{equation}
Hence,
\begin{align}\label{Tstar}
\nonumber \left\|T^*[g] \right\|_p & \le {1 \over 2} \max_\phi \left\| {1 \over \pi} \int_0^\infty  {g(z+r e^{i \phi})-g(z-r e^{i \phi})  \over r } d r  \right\|_p \left| \int_0^\pi e^{-i \phi} d \phi \right| \\
\nonumber &= \max_\phi \left\| {1 \over \pi} \int_0^\infty  {g(z+r e^{i \phi})-g(z-r e^{i \phi})  \over r } d r  \right\|_p \\
&=\max_\phi \left\| \int_{-\infty}^\infty  {g(z+r e^{i \phi}) \over r } d r  \right\|_p=\max_\phi \left\|H[g_\phi] \right\|_p,
\end{align}
where $g_\phi(z)=g(z e^{i \phi})$. The $L_p$-norm in $(\ref{Tstar})$  is two-dimensional, unlike in $(\ref{1DCalZyg})$. However, since
\begin{align}
\nonumber \left\|H[g_\phi] \right\|^p_p&= \int \int_{\field{C}} \left| H[g_\phi](z)\right|^p d x d y=\int_{\infty}^{-\infty}d y  \int_{\infty}^{-\infty} \left| H[g_\phi](x+i y)\right|^p d x  \\
\nonumber  & \le  \int_{\infty}^{-\infty}d y  A^p_p \int_{\infty}^{-\infty} \left| g_\phi(x+i y)\right|^p d x  = A^p_p \int_0^{2 \pi}  \int_0^{\infty} \left| g(\rho e^{i \psi+\phi} ) \right|^p  \rho d \rho d \psi\\
&=A^p_p \int_\phi^{2 \pi+\phi}  \int_0^{\infty} \left| g(\rho e^{i \psi} ) \right|^p  \rho d \rho d \psi=A^p_p \left\| g  \right\|^p_p,
\end{align}
one can get an estimate on the operator $T^*[ \ \cdot \ ]$ in terms of the constant $A_p$ from the one-dimensional bound $(\ref{1DCalZyg})$:
\begin{equation}
\nonumber  \left\|T^*[g]\right\|_p \le  A_p \left\| g\right\|_p \le \cot(\pi/2p) \left\| g \right\|_p.
\end{equation}

Finally,
\begin{equation}
\nonumber \left\|T[g]\right\|_p \le  \cot^2(\pi/2p) \left\| g \right\|_p.
\end{equation}
\end{proof}


\medskip\section{Estimates on the fixed point of the operator ${\bf T_{\bf \mu}}$}\label{Fixed_Point}
\setcounter{equation}{0} 

Theorem $\ref{compactABB}$ suggests the following  approach to the rigorous computer bounds on the solution of the Beltrami equation on the Riemann sphere: Write $\mu$ as $\mu=\kappa+\gamma$ with $\kappa$ and $\gamma$ supported in the neighborhoods of zero and infinity, respectively, find numerically an approximate fixed point of the operator $T_\kappa$ for a $\kappa$ with a large compact support, introduce ``standard sets'' in $L_p(\field{C})$ (we will postpone the discussion of standard sets for later, right now one can think of a standard set in $L_p(\field{C})$ as a ball in $L_p(\field{C})$), choose a standard set that contains our numerical approximation, and show that this standard set contains its image under the operator $T_\kappa$. Then the contraction mapping principle guarantees that the same standard set contains the fixed point $h^*$ of the operator $T_\kappa$. Finally, compute $P_\kappa[h^*]+id$, find bounds on this operation, and estimate the error due to ignoring the differential $\gamma$. 

We will start with a Beltrami differential $\mu \in L_\infty(\hat{\field{C}})$, $\|\mu \|_\infty \le K <1$, given by $\mu=\kappa+ \gamma$ where $\kappa$ and $\gamma$ are compactly supported in a large disk $\field{D}_R$  and a small neighborhood of infinity $\hat{\field{C}} \setminus \field{D}_R$, respectively.   We will further assume that $\kappa$ consists of a ``known'' part $\nu$ (for example, representable as a finite Fourier series, see Appendix), and an ``error''  $\eta$ with a bounded $L_p$ norm: $\kappa=\nu+\eta$, $\left\| \eta \right\|_{p}< \delta$ for some sufficiently small $\delta$.

Let  $h^* \in L_p(\field{C})$ be an approximate fixed point of the operator  $T_\nu$ on $L_p(\field{C})$. Also, let $B_p(h^*,\epsilon)$ denote a ball around $h^*$ in $L_p(\field{C})$ of radius $\epsilon$, i.e the set of all functions $h=h^*+f$ in $L_p(\field{C})$ such that $\left\|f\right\|_p < \epsilon$. By contraction mapping principle the ball $B_p(h^*,\epsilon)$ contains a fixed point of $T_{\nu+\eta}$ if $T_{\nu+\eta}[B_p(h^*,\epsilon)] \subset B_p(h^*,\epsilon)$. The following estimate is straightforward:
\begin{align}\label{Testimate_1}
\nonumber \left\| T_{\nu +\eta} (h)-T_\nu(h^*) \right\|_p&\!=\left\| T [\eta (h^*\!+1) +  (\nu +\eta) f] \right\|_p \le C_p \left\|\eta (h^*\!+1) + (\nu+\eta) f \right\|_p \\
 \nonumber & \le C_p \! \left[ ( \esssup_{\field{D}_R} |h^*\!+1|^p)^{1 \over p}  \left\|\eta\right\|_{p} + (\esssup_{\field{D}_R} |\nu+\eta|^p)^{1 \over p} \left\|f \right\|_p \right]\\
 & \le C_p \! \left[  \left\|\eta\right\|_{p} \sup_{\field{D}_R} |h^*+1| + (\esssup_{\field{D}_R} |\nu\!+\!\eta|^p)^{1 \over p} \left\|f \right\|_p \right],\
\end{align}
where we have used the H\"older inequality and the Calderon-Zygmund lemma. One can use the fact that for a fixed $p<\infty$ and any measurable set $E$ with a finite measure
\begin{equation}
\nonumber \esssup_E |f|^p =\lim_{n \rightarrow \infty} \left[ \int_E |f|^{p n} \right]^{1/n}=\lim_{n \rightarrow \infty} \| f \|^p_{E,p n}= \left( \esssup_E |f|\right)^p,
\end{equation}
to rewrite $(\ref{Testimate_1})$ as
\begin{equation}\label{Testimate}
\left\| T_{\nu +\eta} (h)-T_\nu(h^*) \right\|_p= C_p \left[ \delta \ \sup_{\field{D}_R} {|h^* +1|}+K  \epsilon, \right],
\end{equation}

The above estimate is just one of the several possible ways to bound the difference  $T_{\nu+\eta} (h)-T_\nu(h^*)$ in terms of smallness of $\eta$ and $f$. The choice of norms that measure this smallness, in particular $\| \eta \|_p$ (as opposed to a seemingly more natural $\|\eta\|_\infty$, is specifically adopted to our model uniformization problem (see  Section $\ref{modelproblem}$)).

The estimate $(\ref{Testimate})$ gives a bound on the image of a standard set under $T_{\nu+\eta}$: If one can show that the ball $B_p(T_\nu(h^*),C_p \epsilon')$, $\epsilon'= \delta \  \sup_{\field{D}_R} {|h^* +1|}+K  \epsilon$, is contained in $B_p(h^*,\epsilon)$ for some $p$ such that $K C_p <1$, then the fixed point of the operator $T_{\nu+\eta}$ lies in $B_p(h^*,\epsilon)$.

\medskip\section{Estimates on the solution of the Beltrami equation on a compact} \label{Compact}
\setcounter{equation}{0} 

Define
\begin{equation}
\nonumber g_*^{\nu}(z)=z+P_\nu [h^*](z).
\end{equation}

The proximity of the exact solution of the Beltrami equation on the compact to $g^\nu_*(z)$ is described in the lemma below. 

\begin{lemma}\label{lemmacompact}
  Let the Beltrami differential $\nu+\eta$ be compactly supported in $\field{D}_R$ and  satisfy $\left\|\nu+\eta\right\|_\infty \le K < 1$. Also, given an integer $p>2$, such that $K C_p <1$, let $\left\| \eta \right\|_p < \delta$. Furthermore, assume that the fixed point $h$ of the operator $T_{\nu+\eta}$ lies in  an $\epsilon$-neighborhood of $h^* \in L_p(\field{C})$, i.e. $h \in B_p(h^*,\epsilon)$. Then
\begin{equation}\label{gestimate}
\left|g^{\nu+\eta}(z)-g^\nu_*(z)\right| \le A \left[ \delta  \sup_{\field{D}_R} \left| h^* +1\right| + \epsilon K\right] |z|^{1-2/p},
\end{equation}
where
\begin{equation}\label{constA}
A={1 \over \pi^{1/p}} \left[ 4 { 3 ^{2 q-2} \over 2^q (2-q)} +{25 \over 36} 3^{2 q}-23^{2q-2}+ {\left( {4 \over 9}\right)^{1-q} \over q-1 } \right]^{1 / q},
\end{equation}
and $q$ is the conjugate exponent of $p$.
\end{lemma}
\begin{proof}
Recall that the solution of the Beltrami equation is given by formula $(\ref{solution})$. A rather standard method (see, for instance \cite{Carleson} ) enables ones to get easy H\"older estimates on $\left|g^{\nu+\eta}(z)-g^\nu_*(z)\right|=\left|P_{\nu  + \eta}[h](z) -  P_\nu [h^*](z)\right|$. First,
\begin{align}
\nonumber  \left|g^{\nu+\eta}(z)-g^\nu_*(z)\right|&=\left|P [\eta (h^* + 1) +  (\nu +  \eta) (h-h^*)](z)\right|\\
\nonumber &= {1 \over 2 \pi} \left| \int \int_{\field{D}_R} {\left[ \eta(\xi) (h^*(\xi)  + 1) +  (\nu(\xi)  + \eta(\xi)) f(\xi ) \right] z \over (\xi  -  z) \xi }  d \bar{\xi}  \wedge  d \xi \right| \\
\nonumber & \le {1 \over \pi } \left\|\eta (h^*  + 1) +  (\nu  + \eta) (h-h^*)\right\|_p \! \left[ \int \!\!\int_{\field{D}_R}\!\!{ |z|^q r d r d \phi  \over |(\xi \! -\!  z) \xi|^q }  \right]^{1/q}\!\!\!\!,
\end{align}
where $q$ is the conjugate exponent of $p$ and $f=h-h^*$. Next, split the last integral into three pieces: $D_1=\left\{ \xi \in \field{C}:|\xi-z|<|z|/a \right\}$, $D_2=\left\{ \xi \in \field{C}: |\xi|<|z|/a \right\}$, for some $a \ge 2$,  and $D_3=\field{D}_R \setminus D_1 \cup D_2 $. The first integral over $D_1$, where $|\xi|> (a-1)|z|/a$, can be estimated in the following way:
\begin{equation}
 \nonumber \int  \!\! \! \int_{D_1} \left| {z \over (\xi -z) \xi} \right|^q  r d r d \phi  \le \left({a \over a-1}\right)^q \!\! \int \!\!\! \int_{D_1} { 1 \over |\xi -z|^q } r d r d \phi  \le {  a^{2q-2} \over (a-1)^q} {2 \pi \over 2-q} |z|^{2-q}.
\end{equation}

The integral over $D_2$ is controlled in a similar way. Next, let $D'_3=\left\{ \xi \in \field{C}\right. :|\xi-z/2|\le \left.b |z| \right\} \setminus D_1 \cup D_2 $, for some $b\ge 1/2+1/a$, and $D''_3=\left\{ \xi \in \field{C}\right. :|\xi-z/2|>b |z|, \ |\xi-z/2|<\left. R+|z|/2\right\}$, then the integral over $D_3$ can be  bounded as follows:
\begin{align}
\nonumber \int  \!\!\!\!  \int_{D_3} \left| {z \over (\xi -z) \xi} \right|^q  r d r d \phi &\le \int \!\!\!\!\int_{D'_3}  {|z|^q r d r d \phi  \over |(\xi -z)\xi|^q } + \!\int \!\!\!\! \int_{D''_3} {|z|^q r d r d \phi  \over |(\xi -{z \over 2}) (\xi+{z \over 2})|^q} \\
\nonumber &\le \int \!\!\!\!\int_{D'_3}  {|z|^q r d r d \phi  \over \left|{z \over a }\right|^q \left|{z \over a }\right|^q } + 2 \pi \int^{R+{|z| \over 2}}_{b |z|} {|z|^q r d r  \over (r^2 -{|z|^2 \over 4})^q} \\
 \nonumber &\le a^{2 q} \pi \left(b^2- {2 \over a^2} \right) |z|^{2-q}+ \pi { \left(b^2- {1 \over 4}\right)^{1-q} \over q-1} |z|^{2-q}.
\end{align}

Therefore,
\begin{equation}
\nonumber \left|P [\eta (h^*  +  1)  + (\nu  +  \eta) (h-h^*)](z)\right|  \le A  \left\|\eta (h^* + 1) +  (\nu  + \eta) (h-h^*) \right\|_{p} |z|^{1-{2 \over p}}, \\
\end{equation}
where
\begin{align}\label{bound1}
  \nonumber A=&{1 \over \pi^{1/p} } \left[{  a^{2q-2} \over (a-1)^q} {4 \over 2-q} + a^{2 q} b^2-2 a^{2 q-2} + { \left(b^2- {1 \over 4}\right)^{1-q} \over q-1 }   \right]^{1 \over q}.
\end{align}

The claim of the Lemma follows after one makes a convenient choice of $a=3$, $b=1/2+1/a$. For this choice, $A$ is close to its minimum for small values of $p-2$.
\end{proof}

We recap our procedure upto this point: one first finds a ball in $L_p(\field{C})$ around the approximate solution $h^*$ which is mapped into itself by the operator $T_\nu$; the exact solution $g^{\nu+\eta}$ of the Beltrami equation $g^{\nu+\eta}_{\bar{z}}=(\nu+\eta) g^{\nu+\eta}_z$ satisfies the bound from the above Lemma.

\medskip\section{Estimates on the solution of the Beltrami equation in ${\hat{\field{C}}}$}\label{Plain}
\setcounter{equation}{0} 

Our next step is to estimate the error introduced in ignoring the part of the differential in a small neighborhood of infinity. As explained in the Introduction, define
\begin{equation}\label{beta1}
  \beta=\gamma {g^{\nu+\eta}_z \over \bar{g}^{\nu+\eta}_{\bar{z}}} \circ \left[ g^{\nu+\eta}\right]^{-1}.	
\end{equation}

The solution of the Beltrami equation on a compact is a homeomorphism, hence the continuous inverse of $g^{\nu+\eta}$ is well defined. The differential $\beta$ is again supported in a neighborhood of infinity. Here, rather then getting estimates on the size of this neighborhood from general considerations, we can use our knowledge of  $g^\nu_*$, together with estimate $(\ref{gestimate})$, to obtain a bound on $r$ such that $g^{\nu+\eta}(\field{D}_R)$ contains $\field{D}_{r}$. This amounts to evaluating  $g^\nu_*$ on any collection of standard sets (disks) $Z_k$ in $\field{C}$ that cover the circle $\left\{z \in \field{C}: |z|=R\right\}$ and choosing the infimum of $\left| g^\nu_*(Z_k) \right|$. 

It follows immediately that the differential $\tilde{\beta}(z)=\beta(R^2/z)z^2/\bar{z}^2$ where $\beta$ is defined in $(\ref{beta1})$ is compactly supported in a disk of radius $R^2/r$. In the following Lemma and in the rest of the paper, $\left\| \cdot  \right\|_{E,p}$ will signify the $L_p$-norm over a measurable set $E$.
\begin{lemma} \label{gtilda}
Let $\gamma \in L_\infty(\hat{\field{C}})$ be compactly supported in $\hat{\field{C}} \setminus \field{D}_R$, and satisfy $\| \gamma \|_\infty \le K<1$. Furthermore, let $g^{\tilde{\beta}}$ be the solution of the Beltrami equation $g^{\tilde{\beta}}_{\bar{z}}=\tilde{\beta}g^{\tilde{\beta}}_z$ where $\tilde{\beta}$, given by $(\ref{beta1})$ and $(\ref{betatilde})$, is compactly supported in $\field{D}_{R^2/r}$. Then for all $p>2$ such that $K C_p <1$
\begin{equation}
\nonumber |g^{\tilde{\beta}}(z)-z| \le C |z|^{1-2/p},
\end{equation}
where 
\begin{equation}\label{constC}
C=\pi^{1/p} A K { R^{4/p} ( 2-K C_p ) \over  r^{2/p}( 1-K C_p)}.
\end{equation}
\end{lemma}
\begin{proof}
Notice, that 
\begin{equation}
\nonumber \| \tilde{\beta} \|_\infty=\left\| \beta \right\|_\infty=\left\| \gamma { g^{\nu+\eta}_z \over \bar{g}^{\nu+\eta}_{\bar{z}} }\right\|_\infty=\left\| \gamma \right\|_\infty \le K.
\end{equation}
By a standard result from the theory of quasiconformal mappings (see, for example \cite{Carleson} ), if $\tilde{h}$ is the fixed point of the operator $T_{\tilde{\beta}}$, then 
\begin{equation}
\nonumber \left\|\tilde{h} \right\|_{\field{D}_{{R^2 \over r},p}} \le { \left( \pi R^4\right)^{1/p} \over r^{2/p} (1-K C_p)}
\end{equation}
for all $p$ for which $K C_p <1$.
By an argument similar to that of Lemma  $\ref{lemmacompact}$
\begin{align}
 \nonumber |g^{\tilde{\beta}}(z)-z| & \le A  \left\| \tilde{\beta} \right\|_\infty \! \left\|\tilde{h}+1\right\|_{\field{D}_{R^2/r},p} |z|^{1-2/p} \le \pi^{1/p} A  K { R^{4/p} ( 2-K C_p ) \over  r^{2/p}( 1-K C_p)}  |z|^{1-{2 \over p}}.
\end{align}
\end{proof}

We can now complete the proof of Theorem \ref{Main_Theorem}. 

\bigskip

\noindent {\it Proof of Theorem \ref{Main_Theorem}.}
Let $\beta$ be as in $(\ref{beta1})$ and let, as before, $\tilde{\beta}(z)=\beta(R^2/z) z^2/\bar{z}^2$. Then the following bound on the solution $g^\beta$ of  $g^\beta_{\bar{z}}=\beta g^\beta_{z}$ is immediate:
\begin{align}
\nonumber |g^\beta(z)-z|&=\left|{R^2 \over g^{\tilde{\beta}}(R^2/z)}-z \right|=\left|{R^2\over z} - g^{\tilde{\beta}}(R^2/z) \right| \left| {z\over g^{\tilde{\beta}}(R^2/z)} \right| \\
 \nonumber & \le C \left|{R^2 \over z} \right|^{1-2/p} {|z| \over \left|{R^2 \over z} \right|- C\left|{R^2 \over z} \right|^{1-2/p}}={C |z|^{1+2/p} \over R^{4/p} -C |z|^{2/p} }
\end{align}
Now, we can demonstrate the claim of the Theorem:
\begin{align}
\nonumber  \left|g^\beta \! \circ \! g^{\nu+\eta}(z) -g^\nu_*(z)\right| & \le  \left|g^{\nu+\eta}(z)-g^\nu_*(z) \right|+ \left|g^{\nu+\eta}(z)-g^\beta \circ g^{\nu+\eta}(z)  \right|\\
\nonumber  & \le  A \epsilon' |z|^{1-2/p}  + {C \left|g^{\nu+\eta}(z)\right|^{1+2/p} \over R^{4/p} -C \left|g^{\nu+\eta}(z)\right|^{2/p}}\\
\nonumber & \le  A \epsilon' |z|^{1-2/p}  + {C \left[ \left|g^\nu_*(z)\right|+A \epsilon' |z|^{1-2/p} \right]^{1+2/p} \over R^{4/p} -C \left[ \left|g^\nu_*(z)\right|+A \epsilon' |z|^{1-2/p} \right]^{2/p}}.
\end{align}
Notice, that the first term is $\epsilon'$-small, while the second term is $1/{r}^{2/p}$-small. 
\begin{flushright}{\qed}\end{flushright}

\medskip   \section{Completion of the fundamental domain} \label{Completion}
\setcounter{equation}{0} 

In this Section we will show how one can obtain bounds on the essential supremum and the $L_p$-norm of the Beltrami differential $(\ref{belt2})$ in the neighborhoods of zero and infinity. 

As we have already mentioned in Section $\ref{modelproblem}$, the parts of the boundaries of the fundamental crescent  adjacent to the indifferent fixed point can be chosen as a piece of an internal ray and its preimage under the branch of $f^{-q_n}$ that fixes $0$. We will now briefly describe the construction of these boundaries and estimate the absolute value of the Beltrami differential for small values of the parameter. 

Let $f$ be an analytic map on $\field{C}$ with a fixed point at $0$ such that $\lambda \overset{\rm def}{=}f'(0)=e^{2 \pi i \theta}$ where $\theta$ is 
of bounded type.  By the Siegel Linearization Theorem, there exists a unique conformal local change of coordinates  $\varphi$,  normalized so that $\varphi(0)=0$ and $\varphi'(0)=1$, that maps the disk $\field{D}_{s}$ of radius $s=|\varphi^{-1}(1)|$ onto the Siegel disk and conjugates  $f$ to a rotation by $\lambda$: $f(\varphi (z))=\varphi(\lambda z)$. By the classical Koebe Distortion Theorem
\begin{equation} 
\nonumber |\varphi'(z)-1| \le {|z/s| \over (1-|z/s|)^3}
\end{equation}  
for all $z \in \field{D}_s$. Therefore, if $z \in \field{D}_\rho$, $\rho < s$, then
\begin{equation}
|\varphi(z)-z|  \le c |z|^2, \quad |\varphi'(z)-1| \le c |z|,
\end{equation}
where $c={1 /  s (1-|\rho / s|)^3}$. 

Choose a sufficiently small $\rho>0$ and define $\varrho$ to be the largest parameter value satisfying $\lambda_n(\varrho) \in \varphi(\partial \field{D}_{\rho})$. Define $d=\varphi^{-1}(\lambda_n (\varrho) )$ and consider a ray passing through the point $\lambda_n (\varrho)$: the image  $\varphi(d t)$ of the short line $d t$, $t \in [0,1]$.
The preimage $f^{-q_n} ( \varphi(d t))$ of the ray passes through the point $f^{-q_n} (\lambda_n(\varrho))$ --- it will serve as a completion of the other boundary curve of the fundamental domain.  As before, it is convenient to parametrize this part of the boundary in terms of the radial coordinate $r$. Here, we will make the choice $t=r^\nu/\varrho^\nu$, the positive power $\nu$ will be specified later. As we have demonstrated, the coordinate $\varphi$ is close to the identity on $\field{D}_{\rho}$, therefore the curve $\tilde{\varphi}$ defined by setting $\tilde{\varphi}(r)=\varphi(d r^\nu /\varrho^\nu)$ satisfies the following set of equations for all $r<\varrho$:
\begin{align}
  \nonumber \tilde{\varphi}(r)&={d \over \varrho^\nu} r^\nu (1+\tilde{\gamma}(r)), \quad  {\rm where} \quad  |\tilde{\gamma}(r)|<c |d|,\\
  \nonumber \tilde{\varphi}'(r)&=\nu {d \over \varrho^\nu} r^{\nu-1}(1+\tilde{\beta}(r)), \quad  {\rm where} \quad  |\tilde{\beta}(r)|< c |d|,\\
\nonumber f^{-q_n} (\tilde{\varphi}(r))&= { \lambda^{-q_n} d \over \varrho^\nu} r^\nu(1+\gamma(r)), \quad  {\rm where} \quad   |\gamma(r)|<c |d|,\\
\nonumber (f^{-q_n} \circ \tilde{\varphi})'(r)& = \nu { \lambda^{-q_n} d \over \varrho^\nu} r^{\nu-1} (1+\beta(r)),\quad  {\rm where} \quad  |\beta(r)|<c |d|,
\end{align}
where, as before, $\lambda=e^{2 \pi i \theta}$.

Next,
\begin{equation}
  \nonumber \tau^{-1}(\tilde{\varphi}(r))={1 \over i a }\left[ \ln{\left({r^\nu \over \varrho^\nu}  -{a_{q_n} \over d (1+\tilde{\gamma}(r))}  \right)}-b-\nu \ln{r \over \varrho}\right]={1 \over i a }\left[\tilde{\sigma}(r) -\nu \ln{r \over \varrho }\right],
\end{equation}
where 
\begin{equation}
\nonumber \tilde{\sigma}(r)=\ln{\left( {r^\nu \over  \varrho^\nu}  -{a_{q_n} \over d (1+\tilde{\gamma}(r))} \right)}-b
\end{equation}
is close to 
$\ln{(-a_{q_n} / d )}-b$ for  $r < \varrho $. Similarly,
\begin{equation}
\nonumber \tau^{-1}(f^{-q_n}(\tilde{\varphi}(r)))={1 \over i a }\left[\sigma(r) -\nu \ln{r\over \varrho}-2 \pi i (-q_n \theta \ {\rm  mod} \  1) \right],
\end{equation}
where 
\begin{equation}
\nonumber \sigma(r)=\ln{\left(\lambda^{-q_n}{ r^\nu \over \varrho^\nu}  - {a_{q_n}  \over d (1 +\gamma(r))} \right)}-b
\end{equation}
is  again close to $\ln{(-a_{q_n}  /d)}-b$. Finally, 
\begin{equation}
\nonumber  g_n(r,\phi)= \!\left[ \eta(-\phi)\!+\!{\phi \over 2 \pi}\right] \! {\sigma(r) -  2 \pi i (-q_n \theta \ {\rm  mod} \  1) \over i a}  + \left[ 1\! -\! \eta(-\phi) \! - \! {\phi \over 2 \pi} \right]\! {\tilde{\sigma}(r) \over i a }  - \nu { \ln{r \over \varrho } \over i a}, 
\end{equation}
for small $r$, and, therefore, its complex dilatation on the complement of the non-negative semi-axis is given  by
\begin{align}\label{mu}
\nonumber \mu(r e^{i\phi})&= e^{2 i \phi} {r \partial_r g_n(r,\phi)+i \partial_{\phi} g_n(r,\phi) \over r 
\partial_r g_n(r,\phi)-i \partial_{\phi} g_n(r,\phi) }\\
&= e^{2 i \phi} {\nu- (-q_n \theta \ {\rm  mod} \  1) -\alpha(r,\phi)- {i \over 2 \pi}  (\sigma(r)-\tilde{\sigma}(r))   \over \nu+(-q_n \theta \ {\rm  mod} \  1) -\alpha(r,\phi) +{i \over 2 \pi} (\sigma(r)-\tilde{\sigma}(r))}, 
\end{align}
where 
\begin{equation}\label{alpha_1}
\alpha(r,\phi)= \left(\eta(-\phi)+{\phi \over 2 \pi}\right) r \sigma'(r)+\left(1 -\eta(-\phi)- {\phi \over 2 \pi} \right) r \tilde{\sigma}'(r).  
\end{equation}

Notice, that if $\nu$ is chosen to be $(-q_n \theta \ {\rm  mod} \  1)$, then the absolute value of this dilatation is close to $0$. To provide a bound on this dilatation for all $r \le \varrho$, one needs to know  functions $r \tilde{\sigma}'(r)$ and $r\sigma'(r)$, as formulas $(\ref{mu})$ and  $(\ref{alpha_1})$ indicate: 
\begin{align}
  \nonumber r \tilde{\sigma}'(r)&=   {d (1 +\tilde{\gamma}(r)) \over d (r^\nu / \varrho^\nu) (1+\tilde{\gamma}(r))-a_{q_n}}\left[\nu {r^\nu \over \varrho^\nu} +\nu {a_{q_n} (\tilde{\beta}(r)-\tilde{\gamma}(r)) \over d  (1+\tilde{\gamma}(r))^2} \right], \\
  \nonumber r \sigma'(r)&= {d (1 +\gamma(r)) \over d \lambda^{-q_n} (r^\nu / \varrho^\nu)(1+\gamma(r))-a_{q_n}}\left[ \lambda^{-q_n} \nu {r^\nu \over \varrho^\nu}+\nu {a_{q_n} (\beta(r)-\gamma(r)) \over d  (1+\gamma(r))^2} \right].
\end{align}

It is not hard to evaluate functions $\tilde{\sigma}$, $\sigma$ and $\alpha$  on a standard set $I \supset (0,\varrho)$, and therefore, obtain a bound on the essential supremum of the Beltrami 
differential $(\ref{mu})$ for $r<\varrho$.

As a next step, we will describe the boundary curves of the fundamental domain for very large parameter values. Let $z=a_{q_n}+w$, 
$\xi=f_z^{q_n}(a_{q_n})$, and let 
$\tilde{f}(w)=f^{q_n}(a_{q_n}+w)-a_{q_n}$.  By the Koenigs Linearization Theorem if 
$|\xi| \ne 0,1$ then there exists a unique local holomorphic change of coordinates $\psi$, given by 
$\psi^{-1}(w)=\lim_{k \mapsto \infty} \xi^k \tilde{f}^{-k}(w)$,
 such that $\psi(0)=0$, $\psi'(0)=1$ and $\tilde{f} (\psi(w))=\psi(\xi z)$.  The latter limit is uniform on any disk $\field{D}_s$  such that  
$|\tilde{f}^{-1}(w)| < \bar{c} |w|$ for all $w \in \field{D}_s$ and some fixed constant $\bar{c}$ satisfying $\bar{c}^2<|\xi|^{-1}<\bar{c}$.

As before, fix a small $\rho<s$, then  according to the Koebe Theorem, there exist a constant $c={1 /  s (1-|\rho / s|)^3}$ such  that 
\begin{equation}
|\psi(w)-w|  \le c |w|^2, \quad |\psi'(w)-1| \le c |w|.
\end{equation}

Define $R$ to be the smallest parameter value satisfying $\lambda_n(R) \in \varphi(\partial \field{D}_{\rho}(a_{q_n}))$, $\varphi(z) \equiv \psi(z-a_{q_n})+a_{q_n}$. Define $d=\varphi^{-1}(\lambda_n (R) )$ and consider the logarithmic spiral  $(d-a_{q_n}) e^{ \ \nu \ln{(r/R)} / (2 \pi)}+a_{q_n}$  passing through the point  $\varphi^{-1}(\lambda_n (R))$ when $r=R$. Here $\nu$ is some complex number, to be specified later, with a negative real part. The image 
\begin{equation}
\nonumber \tilde{\varphi}(r)=\psi\left( (d-a_{q_n})  e^{ {\nu \over 2 \pi}  \ln{r \over R}}\right)+a_{q_n} 
\end{equation}
of this spiral is a simple smooth curve contained in $\varphi(\field{D}_{\rho}(a_{q_n})))$ for all $r>R$.  For the sake of notational compactness below, we will use the same notation as in the case of the indifferent fixed point: of course, the functions and constants  appearing in the two cases are completely unrelated. 
\begin{align}
\nonumber \tilde{\varphi}(r)&=(d-a_{q_n}) e^{{\nu \over  2 \pi} \ln{r \over R} } (1+\tilde{\gamma}(r))+a_{q_n}, \quad  {\rm where} \quad   |\tilde{\gamma}(r)|< c_1 |d-a_{q_n}|, \\
\nonumber \tilde{\varphi}'(r)&=(d-a_{q_n}) {\nu \over 2 \pi r}  e^{{\nu \over 2 \pi} \ln{r \over R}} (1+\tilde{\beta}(r)), \quad  {\rm where} \quad  |\tilde{\beta}(r)|< c_2 |d-a_{q_n}|,
\end{align}
and 
\begin{align}
\nonumber f^{-1} (\tilde{\varphi}(r))&={ d-a_{q_n} \over \xi}  e^{{\nu \over  2 \pi} \ln{r \over R} }(1+\gamma(r))+a_{q_n}, \ {\rm where} \ |\gamma(r)|<c_1 {|d-a_{q_n}| \over |\xi|},\\
\nonumber (f^{-1} \circ \tilde{\varphi})'(r)& =\nu {d-a_{q_n}  \over \xi} {e^{{\nu\over 2 \pi} \ln{r \over R}} \over 2 \pi r}  (1+\beta(r)),\quad  {\rm where} \quad  |\beta(r)|< c_2  {|d-a_{q_n}| \over |\xi|}.
\end{align}

Next,
\begin{align}
\nonumber \tau^{-1} \left( \tilde{\varphi}(r)\right) &= {1 \over i a} \left[ \ln{ \left(1-{a_{q_n} \over (d-a_{q_n}) e^{{\nu \over  2 \pi} \ln{r \over R} } (1+\tilde{\gamma}(r)) +a_{q_n}} \right)} -b \right],\\
\nonumber \tau^{-1}(f^{-q_n} (\tilde{\varphi}(r)) &= {1 \over i a} \left[ \ln{ \left(1-{a_{q_n} \over \xi^{-1} (d-a_{q_n}) e^{{\nu \over  2 \pi} \ln{r \over R} } (1+\gamma(r)) +a_{q_n}} \right)} -b \right],
\end{align}
therefore
\begin{align}
\nonumber  g_n(r,\phi)& =  \tau^{-1} \left( \tilde{\varphi}(r)\right) +  {1 \over i a} \left( \eta(-\phi)+{\phi \over 2 \pi}\right)(\sigma(r)-\ln \xi), 
\end{align}
where
\begin{equation} \label{sigma}
\sigma(r)=\ln\left[{1+\gamma(r) \over 1+\tilde{\gamma}(r)} \cdot  {(d-a_{q_n})(1+\tilde{\gamma}(r)) e^{{ \nu \over 2 \pi}  \ln{r \over R} } +a_{q_n} \over \xi^{-1} (d-a_{q_n})(1+\gamma(r)) e^{{\nu \over 2 \pi}  \ln{r \over R} } +a_{q_n} }  \right].
\end{equation}

 The complex dilatation of $g_n$ outside of the non-negative semi-axis is given by
\begin{equation}\label{belt_3}
 \mu(r e^{i\phi})= e^{2 i \phi}{\nu  \!-\!i \ln{\xi}\! + \! 2 \pi r \tilde{\alpha}(r)\! +\!r \left( 2 \pi  \eta(-\phi)\!+\!\phi\right)  (\alpha(r)\!-\!\tilde{\alpha}(r))\!+\! i \sigma(r) \over \nu \!+\!i \ln{\xi} \!+\!  2 \pi r \tilde{\alpha}(r)\!+\!r \left( 2 \pi \eta(-\phi)\!+\!\phi \right)  (\alpha(r)\!-\!\tilde{\alpha}(r))\!-\!i \sigma(r)}, 
\end{equation}
where $\tilde{\alpha}$ and $\alpha$, the $r$-derivatives of $ i a \tau^{-1}(\tilde{\varphi}(r))$ and $ i a \tau^{-1}(f^{-q_n} (\tilde{\varphi}(r)))$, respectively, are given by 
\begin{align}
 \label{tilde_alpha}\tilde{\alpha}(r)  &={\tilde{\gamma}'(r) \over 1+\tilde{\gamma}(r)}-{(d-a_{q_n})e^{{\nu \over 2 \pi}  \ln{r \over R} } ( \tilde{\gamma}'(r) +(1+ \tilde{\gamma}(r)) {\nu \over 2 \pi r})\over (d-a_{q_n}) (1+\tilde{\gamma}(r))  e^{{\nu \over 2 \pi}  \ln{r \over R} } +a_{q_n}},\\
 \label{alpha_2} \alpha(r)& ={\gamma'(r) \over 1+\gamma(r)}-{\xi^{-1}(d-a_{q_n})e^{{\nu \over 2 \pi}  \ln{r \over R} } ( \gamma'(r) +(1+ \gamma(r)) {\nu \over 2 \pi r})\over \xi^{-1}(d-a_{q_n}) (1+\gamma(r))  e^{{\nu \over 2 \pi}  \ln{r \over R} } +a_{q_n}}.
\end{align}

Formulas $(\ref{sigma})$--$(\ref{alpha_2})$,  allow one to estimate supremum of the Beltrami
differential $(\ref{belt_3})$ for $r>R$. We would like to emphasize that the absolute values of functions $\sigma$, $\alpha$ and  $\tilde{\alpha}$  are small for all $r \ge R$ if $d-a_{q_n}$ is small.

The constant $\nu$ can be now conveniently chosen to be equal to $i \ln{\xi}$.

\medskip   \section{Bounds on the Beltrami differential} \label{Bounds}
\setcounter{equation}{0}

Finally, we would like to provide an estimate on the supremum and the $L_p$ norm of the error of the Beltrami differential for intermediate parameter values: $r \in [\varrho,R]$. 

Write the Beltrami differential $(\ref{belt2})$ as 
\begin{equation}
\mu(r e^{i\phi})=e^{2 i \phi} {\phi \kappa(r)+\zeta(r,\phi)+\xi(r) \over \phi \kappa(r)+\zeta(r,\phi)-\xi(r) }, 
\end{equation}
where 
\begin{align}
  \nonumber \kappa(r)&=r { \lambda'_n (r) \over 2 \pi} \left[  D \tau^{-1} \left( f^{-q_n}(\lambda_n (r)) \right) D f^{-q_n}\left( \lambda_n(r) \right)-D \tau^{-1} \left( \lambda_n (r) \right) \right],\\
\nonumber \xi(r)&={i \over 2 \pi} \tau^{-1}(f^{-q_n}(\lambda_n (r)))- {i \over 2 \pi} \tau^{-1}(\lambda_n (r)),\\
\nonumber \zeta(r,\phi) & =\eta(-\phi)\zeta_-(r)+(1-\eta(-\phi)) \zeta_+(r),\\
\nonumber \zeta_-(r) & =r \lambda'_n(r)   D  \tau^{-1} (f^{-q_n}(\lambda_n(r)))D f^{-q_n}\left( \lambda_n(r) \right), \\
\nonumber \zeta_+(r) &= r \lambda'_n (r) D \tau^{-1} \left( \lambda_n(r) \right).
\end{align}

The Fourier coefficients of this Beltrami differential are given by:
\begin{equation}
\nonumber \mu_k(r)={1 \over 2 \pi} \int_{-\pi}^{\pi} \mu (r e^{i \alpha}) e^{-i k \alpha} d \alpha={1 \over 2 \pi} \int_{-\pi}^{\pi}  {\kappa(r) \alpha +\zeta(r,\alpha)+\xi(r) \over  \kappa(r) \alpha+\zeta(r,\alpha)-\xi(r) } e^{i (2-k) \alpha} d \alpha.
\end{equation}

Therefore, one obtains after integration:
\begin{align}
  \nonumber \mu_2(r)&= 1+{\xi(r) \over \pi \kappa(r)}  \left[ \ln{\left(1+{ \pi \over \alpha_+(r)} \right)}- \ln{\left(1-{\pi \over \alpha_-(r)}\right)} \right], \\
\nonumber \mu_{2-k}(r)&= { e^{ -i k\alpha_+(r) } \xi(r) \over \pi \kappa(r)}  \left[  \Ei{\left(i k \left(\pi +\alpha_+(r) \right) \right)}  -  \Ei{\left(i k \alpha_+(r) \right)} \right]\\
 \nonumber    &+ {  e^{-i k \alpha_-(r) } \xi(r) \over \pi \kappa(r)} \left[\Ei{\left(i k \alpha_-(r)  \right)} - \Ei{\left(i k \left(\alpha_-(r)-\pi \right) \right)} \right]\\
\nonumber &= (-1)^k  {\xi(r) \over \pi \kappa(r)} \left[ \ExpEi{\left(i k \left(\pi +\alpha_+(r) \right) \right)} - \ExpEi{\left(\!i k\left(\alpha_-(r)-\pi \right)  \right)} \right]+\\
\nonumber &+{\xi(r) \over \pi \kappa(r)} \left[  \ExpEi{\left(i k\alpha_-(r)  \right)}-  \ExpEi{\left(i k \alpha_+(r) \right)}\right],
\end{align} 
where $\ExpEi(z)\overset{{\rm def}}=e^{-z}\Ei(z)$, function $\Ei$ being the exponential integral, and
\begin{equation}
\nonumber \alpha_-(r)= {\zeta_-(r)-\xi(r) \over \kappa(r)}, \quad  \alpha_+(r)= {\zeta_+(r)-\xi(r) \over \kappa(r)}.
\end{equation}

A simple inspection shows that $\alpha_+(r)+\pi=\alpha_-(r)-\pi$, therefore 

\begin{equation}\label{mu_2mk}
 \mu_{2-k}(r)= \left\{ 
\begin{array}{cc} 
 1+{\xi(r) \over \pi \kappa(r)}  \left[ \ln{ ( \alpha_-(r) )}-\ln{(\alpha_+(r))}\right], & k=0, \\
\atop {\xi(r) \over \pi \kappa(r)} \left[  \ExpEi{\left(i k\alpha_-(r)  \right)}-  \ExpEi{\left(i k \alpha_+(r) \right)}\right], & k \ne 0. 
\end{array}
\right.
\end{equation}

Function $\ExpEi$ is holomorphic in the complex plane with a branch cut along the positive real axis. Of special interest to us is the fact is that the exponential integral has the following representation at zero:
\begin{equation}\label{zseries}
\Ei(z)=\gamma + \ln{(-z)} + \sum^n_{i=1} {z^i \over  i! i} +E_n(z), \quad |E_n(z)| \le {e^{|z|} |z|^{n+1}  \over (n+1)! (n+1) },
\end{equation}
where $\gamma$ is the Euler constant. Furthermore, function $\ExpEi$ itself admits the following asymptotic series  at infinity (cf \cite{Leb}):
\begin{align}
  \label{asseries} \ExpEi(z)& ={1 \over z} \left[\sum^n_{i=0}{i ! \over z^i}+R_n(z)\right],
\end{align}
where
\begin{equation}\label{remainders}
|R_n(z)|   \le  \left[ \eta(\Re{z})  (|z|-|\Im{z}|)  +  1 \right] {(n+1)! \over |z|^{n+1}}.
\end{equation}

We have used series (\ref{zseries}) and (\ref{asseries}) to provide bounds on the function $\ExpEi$, and, therefore, bounds on the coefficients $\mu_{2-k}$, $k \ne 0$, of the Beltrami differential. Moreover, since the typical values of $\alpha_\pm$ are large, asymptotic series (\ref{asseries}) allows one to compute a bound on the supremum of the Beltrami differential on $\overline{ \field{D}_R \setminus \field{D}_\varrho}$. To establish such bound, first, observe that
\begin{align}\label{mu_2mk}
\nonumber \mu_{2-k} (r)&= {\xi(r) \over \pi \kappa(r)} \left[{1-{i \over  k \alpha_-}-{2 \over k^2 \alpha^2_-}+{3 i \over k^3 \alpha^3_-}+ {4 \over k^4 \alpha^4_-}  + R_5(i k \alpha_-) \over i k \alpha_- (r) } \right. \\
\nonumber & \phantom{={\xi(r) \over \pi \kappa(r)} \left[ \right]  } \left.  -{1-{i \over  k \alpha_+}-{2 \over k^2 \alpha^2_+}+{3 i \over k^3 \alpha^3_+}+ {4 \over k^4 \alpha^4_+}  + R_5(i k \alpha_+) \over i k \alpha_+(r) } \right]\\
\nonumber &= {i \xi(r) \over \pi k } \left[{  1-{i \over  k \alpha_+}-{2 \over k^2 \alpha^2_+}+{3 i \over k^3 \alpha^3_+}+ {4 \over k^4 \alpha^4_+}  + R_5(i k \alpha_+)  \over \zeta_+ (r) -\xi (r) } \right. \\
 & \phantom{={\xi(r) \over \pi \kappa(r) }}  \left. - {  1-{i \over  k \alpha_-}-{2 \over k^2 \alpha^2_-}+{3 i \over k^3 \alpha^3_-}+ {4 \over k^4 \alpha^4_-}  + R_5(i k \alpha_-) \over \zeta_- (r) -\xi(r) }  \right],
\end{align}
for all $k \ne 0$. Therefore, $\mu_{2-k}+\mu_{2+k}=A(r)/k^2+B(r)/k^4 + (R^-(r)+R^+(r))/k^6$ and $\mu_{2-k}-\mu_{2+k}=iC(r)/k+iD(r) /k^3+i E(r)/ k^5+(R^-(r)-R^+(r))/k^6$, where
\begin{align}
\nonumber A(r)&={2 \xi(r) \over \pi} \left({1 \over (\zeta_+(r) -\xi(r)) \alpha_+(r)} - {1 \over (\zeta_-(r) -\xi(r)) \alpha_-(r)}   \right), \\
\nonumber B(r)&={6 \xi(r) \over \pi} \left({1 \over (\zeta_-(r) -\xi(r)) \alpha^3_-(r)} - {1 \over (\zeta_+(r) -\xi(r)) \alpha^3_+(r)}   \right), \\
\nonumber C(r)&={2 \xi(r) \over \pi} \left({1 \over \zeta_+(r) -\xi(r)} - {1 \over \zeta_-(r) -\xi(r)}   \right), \\
\nonumber D(r)&={4 \xi(r) \over \pi} \left({1 \over (\zeta_-(r) -\xi(r)) \alpha^2_-(r)} - {1 \over (\zeta_+(r) -\xi(r)) \alpha^2_+(r)}   \right), \\
\nonumber E(r)&={8 \xi(r) \over \pi} \left({1 \over (\zeta_+(r) -\xi(r)) \alpha^4_+(r)} - {1 \over (\zeta_-(r) -\xi(r)) \alpha^4_-(r)}   \right),
\end{align}
and 
\begin{equation}
\nonumber |R^\pm(r)| \le  \left| {\xi(r) \over \pi} \left({R_5(\mp i \alpha_+) \over \zeta_+(r) -\xi(r)} - { R_5(\mp i \alpha_-) \over \zeta_-(r) -\xi(r)}   \right)\right|.
\end{equation}

Next,
\begin{align}
\nonumber |\mu(z)|  & =  \left| \mu_2(r) e^{2 i \phi} + \sum^\infty_{k=1} \mu_{2+k} (r) e^{ i (2+k) \phi} +\mu_{2-k} (r) e^{ i (2-k) \phi}  \right| \\
\nonumber &\le   \left| \mu_2(r) + \sum^\infty_{k=1} (\mu_{2-k}(r)+\mu_{2+k}(r)) \cos{k \phi} - i (\mu_{2-k}(r)-\mu_{2+k}(r)) \sin{k \phi} \right|\\
\nonumber & \le   \left| \mu_2(r) + \sum^\infty_{k=1} A(r) {\cos{k \phi} \over k^2}+B(r) {\cos{k \phi} \over k^4} + C(r) {\sin{k \phi} \over k} +  D(r) {\sin{k \phi} \over k^3}+ \right. \\
&  \phantom{\le  \mu_2(r) a} \left.  + E(r) {\sin{k \phi} \over k^5}+ R^-(r) {e^{-i k \phi} \over k^6} +R^+(r) {e^{i k \phi}  \over k^6} \right|.
\end{align}

Therefore, $\sup_{z \in \overline{\field{D}_R \setminus \field{D}_\varrho}} |\mu(z)|$ is bounded by the supremum of the following expression on $[\varrho,R] \times (0,2 \pi)$:
\begin{align} \label{mu_sup}
\nonumber  & \left| \mu_{2}(r) + A(r) \left[{\pi^2 \over 6} -{\pi \phi \over 2}+{\phi^2 \over 4} \right] + B(r) \left[ {\pi^4 \over 90 } -{\pi^2 \phi^2 \over 12} +{\pi \phi^3 \over 12} -{\phi^4 \over 48}  \right] +C(r) {\pi - \phi \over 2} \right.\\
\nonumber&  \phantom{ \left| \mu_{2}(r) \right.  } +\left.D(r)\left[{\pi^2 \phi \over 6} -{\pi \phi^2 \over 4}+{\phi^3 \over 12}    \right] +E(r) \left[{\pi^4 \phi  \over 90} -{\pi^2 \phi^3 \over 36}+{\pi \phi^54\over 48}   -{\phi^5 \over 240}  \right] \right| \\
 & \phantom{ \left| \mu_{2}(r) \right.  }+{\pi^6 \over 945} (|R^-(r)|+|R^+(r)|).
\end{align}

Lastly, we will estimate the accuracy of the approximation of the Beltrami differential with a finite Fourier series. Fix some sufficiently large even natural $M$ and represent $\mu$ on $\field{D}_R$ as $\nu+\eta$, where $\nu$ is supported in $\field{D}_R \setminus \field{D}_\varrho$ and $\eta=\mu$ on $\field{D}_\varrho$. Furthermore, we can assume that $\nu$ is representable as a finite Fourier series on $\field{D}_R \setminus \field{D}_\varrho$, and $\eta$ is the ``higher order error''  on $\field{D}_R \setminus \field{D}_\varrho$: 
\begin{align}
\nonumber \nu(r e^{i \phi})&=\sum^{M+2}_{k=2-M}\mu_k(r)e^{i k \phi },\\
\nonumber \eta(r e^{i \phi})&=\sum^{\infty}_{k=M+3}\mu_k(r)e^{i k \phi }+\sum^{1-M}_{k=-\infty}\mu_k(r)e^{i k \phi }=\eta_1(r e^{i \phi})+\eta_2(r e^{i \phi}),\\
\nonumber \eta_1(r e^{i \phi})&=e^{2 i \phi } \sum^\infty_{k=M+1} \left[{A(r) \over k^2}+{B(r) \over k^4} \right] \cos{k \phi} + \left[{C(r)  \over k} +  {D(r) \over k^3}+ {E(r) \over k^5}\right] \sin{k \phi},\\
\nonumber \eta_2(r e^{i \phi})&= e^{2 i \phi } \sum^\infty_{k=M+1} R^-(r) {e^{-i k \phi} \over k^6} +R^+(r) {e^{i k \phi}  \over k^6}.
\end{align}

To obtain a bound on the $L_p$-norm of $\eta$ on $\overline{ \field{D}_R \setminus \field{D}_\varrho }$, one can use the fact that  $\|\eta\|_{\overline{ \field{D}_R \setminus \field{D}_\varrho},p} \le \sup_{z \in \overline{\field{D}_R \setminus \field{D}_\varrho}} |\eta(z)|^{1-2/p}  \| \eta \|^{2/p}_{\overline{ \field{D}_R \setminus \field{D}_\varrho},2}$ for $p>2$, together with the following bound on the squares of the $L_2$ norms of $\eta_1$ and  $\eta_2$:
\begin{align}
\nonumber  \|\eta_1\|^2_2 & \le 2 \pi \int^R_\varrho \sum^{\infty}_{k=M+3}|\mu_k(r)|^2+\sum^{1-M}_{k=-\infty}|\mu_k(r)|^2 r d r\\
\nonumber  &= 2 \pi \int^R_\varrho \sum^{\infty}_{k=M+1}|\mu_{2+k}(r)|^2+ \sum^{\infty}_{k=M+1}|\mu_{2-k}(r)|^2 r d r\\
\nonumber  &= \pi \int^R_\varrho \sum^{\infty}_{k=M+1}|\mu_{2-k}(r)+\mu_{2+k}(r)|^2+|\mu_{2-k}(r)-\mu_{2+k}(r)|^2 r d r\\
\nonumber  &= {\pi \over 2} \sum^{N-1}_{i=0}\left[ |C_i|^2 \zeta_M(2)+\left[ |A_i|^2 +2 \Re{ ( C_i\overline{D_i} )} \right]  \zeta_M(4) \right.  +\\
\nonumber & \phantom{aaaaaaa} \left[|D_i|^2+2 \Re{ (C_i\overline{E_i})} +2 \Re{ (A_i\overline{B_i} )} \right]  \zeta_M(6)+\\
\label{eta1_norm} & \phantom{aaaaaaa}\left. \left[|B_i|^2+2 \Re{ (D_i\overline{E_i} )} \right]  \zeta_M(8)+ |E_i|^2 \zeta_M(10) \right] (\rho^2_{i}-\rho^2_{i+1}),\\
\nonumber \|\eta_2\|^2_2 & \le  2 \pi \int^R_\varrho r d r \ \left[  | R^-(r)|^2 +| R^+(r)|^2 \right] \zeta_M(12)\\
\label{eta2_norm} & =  \pi \zeta_M(12) \sum^{N-1}_{i=0}(\rho^2_{i}-\rho^2_{i+1})  \left[  | R^-_i|^2 +| R^+_i|^2 \right],
\end{align}
where $\zeta_M(n)$ are the $(M+1)$-remainders of the Riemann zeta function: $\zeta_M (n)=\sum^\infty_{i={M+1}} 
i^{-n}$, and $A_i$, $B_i$, $C_i$, $D_i$, $E_i$ and $R^\pm_i$ are bounds on the corresponding functions on intervals 
$(\rho_{i},\rho_{i+1})$ that cover $(\varrho,R)$, $\rho_0=\varrho$ (see Appendix for details). We would like to 
emphasize that the remainder of the zeta function $\zeta_M(2)$ decreases very slowly with a growing $M$. The 
appearance of this function in the leading term of the norm of $\eta$ necessitates working with 
high order truncations of the Fourier series. 

Naturally, bounds (\ref{mu_sup}), (\ref{eta1_norm}) and (\ref{eta2_norm}) can be implemented on a computer.

The $L_p$-norm of $\eta$ on $\field{D}_\varrho$ can be bounded by $\sup_{z \in \field{D}_\varrho}| \mu(z)| \left(\pi \varrho^2 \right)^{1/p}$, where $\sup_{z \in \field{D}_\varrho} |\mu(z)|$ is computable with the help of formula  $(\ref{mu})$.

\medskip   \section{Appendix. Standard sets}\label{AppendixA}
\setcounter{equation}{0} 

Over the course of the last two decades computer-assisted proofs have become a standard mathematical tool in dynamics since they first appeared  in O. E. Lanford's \citeyear{Lan1,Lan2} pioneering works in the early 80-s (also, cf \cite{EKW1}, \cite{EKW2}, \cite{M}, \cite{LR}, \cite{St} and  \cite{Koch}  for some other important computer-assisted results). Here we will give only a very brief description of some ideas: An interested reader is referred to an excellent review \cite{KSW} on the subject. 

Given a set $X$, denote by $\mathcal{P}(X)$ its power set. Let $f$ be a map from a subset $D_f \subset X$ to another set $Y$. We say that $F: D_F \subset \mathcal{P}(X) \mapsto \mathcal{P}(Y)$ is a bound on $f$ if $f(z) \in F(Z)$ whenever $z \in Z \in D_F$. It is customary to construct such bounds within the class of maps $F: \std{(X)} \mapsto \std{(Y)}$, where $\std{(X)}$, the collection of {\it standard sets} in $X$, is a conveniently chosen subset of $\mathcal{P}(X)$.

\subsection{Interval arithmetics in $\field{C}$}
A computer implementation of an arithmetic operation  $r_1 \# r_2$ ($\#$ is $+$, $-$, $*$ or $/$) on two real numbers does not generally yield an exact result. The ``computer'' result is a number representable in a standard IEEE floating point format (cf \cite{IEEE}). Such numbers are commonly referred to as ``representable''. However, it is relatively straightforward to find the smallest representable number larger than $r_1 \# r_2$ and the largest representable number smaller than  $r_1 \# r_2$, if such numbers exist. If either of these numbers, denoted $r_>$ and $r_<$,  does not exist then the pair $(r_1,r_2)$ is considered to be outside of the domain of the map $(r_1,r_2) \mapsto r_1 \# r_2$ and the program is aborted. If such bounds do exist, then one can take the interval $I[r_<,r_>]$ as a bound on the result of the operation. Thus, one is naturally brought to the idea of constructing bounds on functions on real numbers as functions on the set of all representable intervals:
\begin{equation}
\std{(\field{R})}=\{I[x,y] \in \field{R}: x,y - {\rm representable} \}.
\end{equation}    

To obtain a bound on a transcendental function ($\exp$, $\cos$, $\sin$, $\ln$ and others) one can use their Taylor series together with a bound on the remainder. Some algebraic functions, like roots, can be bounded with the help of these transcendental functions.

Arithmetic operations on complex numbers are reduced to those on reals, and so is the evaluation of bounds of functions on $\field{C}$ on  standard sets in $\field{C}$. We  will choose these standard sets to be open disks in the complex plain with representable centers and radii:
\begin{equation}
\std{(\field{C})}=\{B[z,r] \overset{\rm def}=\field{D}_r(z) \subset \field{C}:z=x+i y;  x,y,r - {\rm representable} \}.
\end{equation}

\subsection{Operations on standard sets for analytic functions}
Given a function $f \in {\bf C}_\rho$ we will bound it by a standard set of the following form:
\begin{align}\label{stdA}
\nonumber \!  F\left({\bf B} ,I_g,I_h\right) &=\left\{ f \in {\bf C}_\rho: f(z)=p(z) + z g(z) + z^{N+1} h(z); \phantom{\sum^N_{k=1}} \right.\\
& \left.\phantom{==|} p(z)=\sum^N_{k=1} c_k z^k, c_k \in B_k, \|g\|_\rho \in I_g, \|h\|_\rho \in I_h \right\},
\end{align}
where $N$ is a sufficiently large fixed integer, ${\bf B}$ is the shorthand for the array of $B_k \in \std{(\field{C})}$, $k=1...N$, and $I_g$, $I_h$ are in $\std{(\field{R})}$. The function $g(z)$ is commonly referred to as the ``general error'' (cf \cite{KSW} for the explanation of the necessity of this term), the function  $h(z)$ is usually called  the ``higher order error''.  

Standard sets $(\ref{stdA})$ are well suited for obtaining bounds on addition, multiplication, composition and differentiation of functions in ${\bf C}_\rho$. As an example, we will consider  bounds on multiplication of two analytic functions. Let $f_1$ and $f_2$  be in ${\bf C}_\rho$. Given two bounds on these functions, $F_1 \ni f_1$ and $F_2 \ni f_2$, we would like to obtain a standard set $F \in \std{({\bf C}_\rho)}$ which contains $f_1 \cdot f_2$. Let $p_1=\sum^N_{k=1} a_k z^k$ and $p_2=\sum^N_{k=1} b_k z^k$. Let $p_1(z) \cdot p_2(z)=p (z) + z^{N+1} r(z)$ where the order of $p$ is no more than $N$. Then 
\begin{align}
\nonumber (f_1 \cdot f_2)(z)=p(z)&+z  \left[p_1(z) \cdot g_2(z)+g_1(z) \cdot p_2(z) +z \cdot g_1(z) \cdot g_2(z)  \right] \\
\nonumber & +z^{N+1}  \left( r(z)+z^{N+1} \cdot h_1(z) \cdot h_2(z)+(p_1(z)+zg_1(z)) \cdot h_2(z) \right. \\ 
&\phantom{ +z^{N+1}ar(z)a}  + \left. h_1(z) \cdot (p_2(z)+z g_2(z)) \right).
\end{align}

Here, the terms in the brackets constitute a general error, while those in the parentheses are the higher order error terms. One knows the bounds on all objects that enter these error terms, therefore one can find the interval $I_g$ and $I_h$ that contain the norms of the general and higher order terms of $f_1 \cdot f_2$.

We would like to note that although our programs include all the necessary routines that supply bounds on operations with standard sets $(\ref{stdA})$, these routines have not been fully used in the proof of the Main Example. This Example serving only demonstrative purposes, the analytic function appearing in it is a quadratic polynomial. Operations with this polynomial do not require estimates on general and higher order terms (if $N$ is chosen to be sufficiently high).

\subsection{Standard sets and operators in ${\bf L_{\bf p}}(\field{C})$}
We have adopted the strategy of approximating a function $h(r e^{i \phi})=\sum^{\infty}_{k=-\infty}h_k(r)e^{i k \phi}$  in $L_p(\field{C})$ by a finite Fourier series with piecewise-constant compactly supported coefficients. Given a collection of real numbers  $0=r_0 < r_1 < \ldots < r_{N-1} <  r_N < r_{N+1} < \ldots < r_M$ that specify a radial grid, we will denote the midpoint of the interval $[r_m,r_{m+1}]$ by $\varrho_m$; $\varrho_{-1}$ will be $0$ by definition. We will make the following choice of the standard sets in $L_p(\field{C})$ associated with a grid $\{r_m\}^{M}_{m=0}$.
\begin{align}\label{stdLp}
  \nonumber  H \left({\bf S},I_e\right) & = \left\{  h=s+e \in L_p(\field{C}) :   s(r e^{i \phi})=  \!\!\! \sum^{M_2}_{k=-M_1}  \!\!\! s_k(r) e^{i k \phi},   s_k(r) \in  S_{k,m} \ {\rm for } \ r \in  \right.\\
 & \left.\phantom{ \sum^M_{k=-M}\!\!\!\!\!\!\!} [\varrho_{m-1},\varrho_m), m \le N, \ {\rm and} \  s_k(r)=0 \ {\rm for}  \ r \ge \varrho_N;  \ \|e\|_p \in I_e \right\}.
\end{align}
Here, $M_1$ and $M_2$ are fixed positive integers, and ${\bf S}$ is the shorthand for the array of $S_{k,m} \in \std{(\field{C})}$.

Let $h=s+e \in L_p(\field{C})$ lie in a standard set  $H\left({\bf S},I_e\right)$. We will now describe the standard set  that contain $T[h]$ and $P[h]$. With this goal in mind, we will describe the action of the Hilbert and Cauchy operators on  piecewise-constant functions.

 Let $s$ be a compactly supported function in $L_p(\field{C})$ given by a finite Fourier series with piecewise-constant coefficients:
\begin{equation} \label{lseries}
 s(r e^{i \phi})=\sum^{M_2}_{k=-M_1} s_k(r)e^{i k \phi} , \ s_k(r)= \left\{  s_{k,m}, \ r \in [\varrho_{m-1},\varrho_m),\ 0 \le m \le N, \atop 0, \  r \ge  \varrho_N, \right.
\end{equation}

The action of the Hilbert operator on such function again can be represented as a Fourier series:
\begin{equation}
T[s](r e^{i \phi})=\sum^{M_2}_{k=-M_1} c_k(r) e^{i k \phi}.
\end{equation}

The coefficients of this series are given  by the following set of expressions (cf \cite{GaiKhmel} and \cite{Daripa3}):
\begin{itemize}
\item[1)]$r=0$
\begin{equation}
\nonumber c_k(0)=\left\{ 
\begin{array}{cc} 
  c_0(r_1)-s_2(r_1) -2 s_2(r_1) \ln{r_1 \over \varrho_0 }, &  k=0, \\
   0,&  k \ne 0.
\end{array}
  \right.
\end{equation}
\item[2)] $r \in [\varrho_{m-1},\varrho_m)$, $0 \le m \le N$,
\begin{equation}
\nonumber c_k(r)=\!\left\{ 
\begin{array}{cc} 
\left( {r \over r_{m+1}} \right)^k \left(c_k(r_{m+1})-s_{k+2}(r_{m+1})\right)-d_{k,m}(r) +s_{k+2}(r), &  k \ge 1,\\
 \! c_0(r_{m+1})-s_2(r_{m+1}) -2 s_2(r_{m+1}) \ln{r_{m+1} \over \varrho_m }-2 s_2(r_m) \ln{\varrho_m \over r }+s_2(r), &  k=0, \\
\left( {r \over r_{m-1}} \right)^k \left(c_k(r_{m-1})-s_{k+2}(r_{m-1})\right)+b_{k,m}(r) +s_{k+2}(r),&  k < 0,
\end{array}
  \right. 
\end{equation}
where 
\begin{align}
  \nonumber d_{k,m}(r)&=2 {k+1 \over k} \left[s_{k+2}(r_{m+1}) \left({r^k \over \varrho^k_m}-{r^k \over r^k_{m+1}} \ \right) + s_{k+2}(r_m) \left(1-{r^k \over \varrho^k_m } \right) \right], \\
  \nonumber b_{k,m}(r)&=2 {k+1 \over k} \left[s_{k+2}(r_{m-1}) \left({r^k \over r^k_{m-1}}-{r^k \over \varrho^k_{m-1}} \right) + s_{k+2}(r_m) \left({r^k \over \varrho^k_{m-1} }-1\right) \right].
\end{align}

\item[3)] $r \in [\varrho_{m-1},\varrho_m)$, $N < m \le M$
  \begin{equation}
\nonumber c_k(r)=\left\{ 
\begin{array}{cc} 
  0,&  k \ne 0,\\
 {r^k \over r^k_{m-1}} c_k (r_{m-1}), &  k < 0.
\end{array}
  \right. 
  \end{equation}
\end{itemize}

To find the standard set $H({\bf C},I[0,0])$ that contains $T[s]$ one has to identify  the standard sets $C_{k,m}$ such that $c_k(r)\in  C_{k,m}$ whenever   $r \in [\varrho_{m-1},\varrho_m)$. This amounts to evaluation of $c_k(r)$ on an interval $I_m \in \std{(\field{R})}$  containing  $[\varrho_{m-1},\varrho_m)$.

Finally, the interval $I_{T[e]}$ that contains $\|T[e]\|_p$ can be readily found with the help the Calderon-Zygmund lemma: It is the unique standard set with the smallest length in $\std{(\field{R})}$ that contains all real numbers $C_p x$ for all $x \in I_e$.

In a similar fashion the action of the Cauchy transform on a piecewise-constant function $(\ref{lseries})$ is given by the following set of equations:

\begin{itemize}
\item[1)]$r=0$
\begin{equation}	
\nonumber c_k(0)=\left\{		 
\begin{array}{cc} 		
  c_0(r_1)-2 s_1(r_1) \varrho_0 -2 s_1(r_1) (r_m-\varrho_0), &  k=0, \\
   0,&  k \ne 0.
\end{array}	
  \right.
\end{equation}
\item[2)] $r \in [\varrho_{m-1},\varrho_m)$, $0 \le m \le N$,
\begin{equation}
\nonumber c_k(r)=\left\{ 
\begin{array}{cc} 
\left( {r \over r_{m+1}} \right)^k c_k(r_{m+1})-d_{k,m}(r), &  k > 1,\\
 {r \over r_{m+1}} c_1(r_{m+1}) -2 r s_2(r_m) \ln{\varrho_m \over r}-2 r s_2(r_{m+1})\ln{r_{m+1} \over \varrho_m} , &  k=1, \\
\left( {r \over r_{m-1}} \right)^k c_k(r_{m-1})+b_{k,m}(r),&  k <= 0,
\end{array}
  \right. 
\end{equation}
where 
\begin{align}
  \nonumber d_{k,m}(r)&={2 r \over 1-k} \left[s_{k+1}(r_m) \left({r^{k-1} \over \varrho^{k-1}_m}-1 \right) + s_{k+1}(r_{m+1}) \left({r^{k-1} \over r^{k-1}_{m+1} }-{r^{k-1} \over \varrho^{k-1}_m } \right) \right], \\
  \nonumber b_{k,m}(r)&={2 r \over 1-k} \left[s_{k+1}(r_{m-1}) \left({\varrho^{1-k}_{m-1} \over r^{1-k}}-{r^{1-k}_{m-1} \over r^{1-k}} \right) + s_{k+1}(r_m) \left(1-{\varrho^{1-k}_{m-1} \over r^{1-k}}\right) \right].
\end{align}

\item[3)] $r \in [\varrho_{m-1},\varrho_m)$, $N < m \le M$
  \begin{equation}
\nonumber c_k(r)=\left\{ 
\begin{array}{cc} 
  0,&  k \ne 0,\\
 {r^k \over r^k_{m-1}} c_k (r_{m-1}), &  k < 0.
\end{array}
  \right. 
  \end{equation}
\end{itemize}

The bound on $P[s]$ is found similarly to the analogous estimate for the Hilbert transform. Notice, that we do not need the bound on $P[e]$: According to Theorem \ref{Main_Theorem}, one only requires bounds on the action of the Cauchy transform on the approximate fixed point $h^*$ of the iteration scheme. Such approximate fixed point is always given as a finite Fourier series.

\medskip   
\section{Acknowledgments.} 

The author would like to thank the Department of Mathematics of the University of Texas at Austin for making their AMD64 computational cluster available for this project. 

Many thanks go out to Michael Yampolsky for his numerous helpful comments and his ideas concerning the improvement of the computer implementation of cylinder renormalization.

The author grieves the untimely passage of Dima Khmelev. Dima was an important member of the cylinder renormalization project that necessitated the present work and contributed to our studies of the numerics of the Beltrami equation in a crucial way.

\end{document}